\RequirePackage{etoolbox}
\csdef{input@path}{{style/}{graphics/}}
\documentclass[MSNbibl,number,citesort,seceqn,dvips]{arxbj}
\usepackage{mathbh}
\usepackage{upgreek}

\aid{0}
\volume{21}
\issue{1}
\pubyear{2015}
\firstpage{176}
\lastpage{208}
\doi{10.3150/13-BEJ564} 

\makeatletter
\newcommand{\rrvert}{\vert}
\newcommand{\llvert}{\vert}
\newtheorem{them}{Theorem}
\newtheorem{cor}{Corollary}
\newtheorem{lemma}{Lemma}[section]
\newproclaim{MarginAssumption}{Margin assumption}
\newproclaim{NoiseAssumption}{Noise assumption}
\newproclaim{PluginAssumption}{Plug-in assumption}
\newproclaim{Definition}{Definition}
\newproclaim{KernelAssumption}{Kernel assumption}

\renewcommand{\pi}{\uppi}
\renewcommand{\epsilon}{\varepsilon}
\newcommand{\rodsigma}{\scriptsize{\overrightarrow{\sigma\vphantom{A}}}}
\newcommand{\rodsigmas}{{\fontsize{6}{8}\selectfont{\overrightarrow{\sigma\vphantom{A}}}}}
\newcommand{\tiksigma}{\overrightarrow{\sigma\vphantom{A}}}

\newcommand{\eqref}[1]{(\ref{#1})}
\newcommand{\trup}[2]{{#1}/{#2}}
\newcommand{\GG}{\mathcal{G}}
\newcommand{\FF}{\mathcal{F}}
\newcommand{\no}{\Vert}
\newcommand{\ind}{\mathbh{1}}
\newcommand{\R}{\mathbb{R}}
\renewcommand{\P}{\mathbb{P}}
\newcommand{\E}{\mathbb{E}}
\makeatother

\begin{document}
\begin{frontmatter}

\title{Minimax fast rates for discriminant analysis with errors in variables}
\runtitle{Minimax rates for noisy discriminant analysis}

\begin{aug}
\author[1]{\inits{S.}\fnms{S\'ebastien} \snm{Loustau}\corref{}\thanksref{1}\ead[label=e1]{loustau@math.univ-angers.fr}} \and
\author[2]{\inits{C.}\fnms{Cl\'ement} \snm{Marteau}\thanksref{2}\ead[label=e2]{marteau@insa-toulouse.fr}}
\address[1]{Laboratoire Angevin de Recherche en Maths, Universit\'e
d'Angers, 2 Boulevard Lavoisier, 49045 Angers Cedex 01, France. \printead{e1}}
\address[2]{Institut de Math\'{e}matiques de Toulouse, INSA de Toulouse,
Universit\'{e} de Toulouse, 118 route de Narbonne, F-31062 Toulouse
Cedex 4, France. \printead{e2}}
\dedicated{This work is dedicated to the memory of our mentor and colleague Laurent Cavalier, who suddenly passed away on January 2014}
\end{aug}

\received{\smonth{6} \syear{2012}}
\revised{\smonth{5} \syear{2013}}

%
\begin{abstract}
The effect of measurement errors in discriminant analysis
is investigated. Given observations $Z=X+\varepsilon$, where $\varepsilon$
denotes a random noise, the goal is to predict the density of $X$
among two possible candidates $f$ and $g$. We suppose that we have
at our disposal two learning samples. The aim is to approach the best
possible decision rule $G^\star$ defined as a minimizer of the Bayes
risk.

In the free-noise case ($\varepsilon=0$), minimax fast rates of
convergence are well-known under the margin assumption in discriminant
analysis (see (\textit{Ann. Statist.}
\textbf{27} (1999)
1808--1829)) or in the more general classification
framework (see (\textit{Ann. Statist.}
\textbf{35} (2002)
608--633, \textit{Ann. Statist.}
\textbf{32} (2004)
135--166)). In this paper, we intend to
establish similar results in the noisy case, that is, when dealing
with errors in variables. We prove minimax lower bounds for this
problem and explain how can these rates be attained, using in
particular an Empirical Risk Minimizer (ERM) method based on
deconvolution kernel estimators.
\end{abstract}

%

\end{frontmatter}

\section{Introduction}

In the problem of discriminant analysis, we usually observe two i.i.d.
samples $X_1^{(1)},\ldots, X_n^{(1)}$ and $X_1^{(2)},\ldots, X_m^{(2)}$. Each
observation $X_j^{(i)}\in\R^d$ is assumed to admit a density with
respect to\vspace*{2pt} a \mbox{$\sigma$-finite} measure $Q$, dominated by the Lebesgue
measure. This density will be denoted by $f$ if the observation belongs
to the first set (i.e., when $i=1$) or $g$ in the other case. Our aim
is to infer the density of a new incoming observation $X$. This problem
can be considered as a particular case of the more general and
extensively studied binary classification problem (see \cite{jaune}
for a detailed introduction or \cite{survey} for a concise survey).

In this framework, a decision rule or classifier can be identified with
a set $G \subset\R^d$, which attributes $X$ to $f$ if $X\in G$ and to
$g$ otherwise. Then, we can associate to each classifier $G$ its
corresponding Bayes risk $R_K(G)$ defined as:
%
\begin{equation}
\label{risk} R_K(G)=\frac{1}{2} \biggl[\int
_{K/G}f(x)\,\mathrm{d}Q(x)+\int_Gg(x)
\,\mathrm{d}Q(x) \biggr],
\end{equation}
where we restrict the problem to a compact set $K\subset\R^d$. The
minimizer of the Bayes risk (the best possible classifier for this
criterion) is given by:
%
\begin{equation}
\label{bayes} G_K^\star=\bigl\{x\in K\dvtx f(x)\geq g(x)
\bigr\}, 
\end{equation}
where the infimum is taken over all subsets of $K$. The Bayes
classifier is obviously unknown since it explicitly depends on the
couple $(f,g)$. The goal is thus to estimate $G_K^\star$ thanks to a
classifier $\hat{G}_{n,m}$ based on the two learning samples.

The risk minimizer \eqref{bayes} has attracted many attentions in the
last two decades because it involves a quantity of applied motivating
examples, including pattern recognition, spam filtering, or medical
diagnostic. However, in many real-world problems, direct observations
are not available and measurement errors occur. As a result, it could
be interesting to take into account this problem into the
classification task. In this paper, we propose to estimate the Bayes
classifier $G_K^\star$ defined in (\ref{bayes}) thanks to noisy
samples. For all $i\in\{1,2 \}$, we assume that we observe:
%
%
\begin{equation}
\label{eq:modelb} Z_j^{(i)} = X_j^{(i)} +
\epsilon_j^{(i)},\qquad j=1,\ldots, n_i,
\end{equation}
instead of the $X_j^{(i)}$, where in the sequel $n_1=n$ and $n_2=m$.
The $\epsilon_j^{(i)}$ denotes i.i.d. random variables expressing
measurement errors. We will see in this work that we are facing an
inverse problem, and more precisely a deconvolution problem. Indeed,
assume that for all $x\in\mathbb{R}^d$, $\mathrm{d}Q(x)=\mu(x) \,\mathrm{d}x$ for some
bounded function $\mu$. If $\epsilon$ admits a density $\eta$ with
respect to the Lebesgue measure, then the corresponding density of the
$Z_j^{(i)}$ is the convolution product $(f\cdot\mu)*\eta$ if $i=1$ or
$(g\cdot\mu)*\eta$ if $i=2$. This property gives rise to a deconvolution
step in the estimation procedure. Deconvolution problems arise in many
fields where data are obtained with measurement errors and are at the
core of several nonparametric statistical studies. For a general review
of the possible methodologies associated to these problems, we may
mention for instance \cite{meister}. More specifically, we refer to
\cite{fan} in density estimation, \cite{Carroll_Delaigle_Hall} for
nonparametric prediction or \cite{Butucea} where goodness-of-fit
tests are constructed in the presence of noise. The main key of all
these studies is to construct a deconvolution kernel which may allow to
annihilate the noise $\epsilon$. More details on the construction of
such objects are provided in Section~\ref{section:upperbounds}. It is
important to note that in this discriminant analysis setup, or more
generally in classification, there is up to our knowledge no such a
work. The aim of this article is to describe minimax rates of
convergence in noisy discriminant analysis under the Margin
assumption.

In the free-noise case, that is, when $\epsilon=0$, \cite{mammen} has
attracted the attention on minimax fast rates of convergence (i.e.,
faster than $n^{-\trup{1}{2}}$). In particular, they propose a
classifier $\hat G_{n,m}$ satisfying
%
\begin{equation}
\label{fastrates} \sup_{G_K^\star\in\mathcal{G}(\alpha,\rho)} \E \bigl[ R_K(\hat
{G}_{n,m})-R_K\bigl(G_K^\star
\bigr) \bigr] \leq C (n\wedge m)^{-\trup{(\alpha
+1)}{(2+\alpha+\rho\alpha)}},
\end{equation}
for some positive constant $C$. Here, $\mathcal{G}(\alpha,\rho)$
denotes a nonparametric set of candidates $G_K^\star$ with complexity
$\rho>0$ and margin parameter $\alpha\geq0$ (see Section~\ref{sec2.1} for a
precise definition). In \eqref{fastrates}, the complexity parameter
$\rho>0$ is related to the notion of entropy with bracketing whereas
the margin is used to relate the variance to the expectation. It allows
\cite{mammen} to get improved bounds using the so-called peeling
technique of \cite{vdg}. This result is at the origin of a recent and
vast literature on fast rates of convergence in classification (see,
for instance, \cite{nedelec,AT}) or in general statistical learning
(see \cite{kolt}). In these\vadjust{\goodbreak} papers, the complexity assumption can be
of two forms: a~geometric assumption over the class of candidates
$G_K^\star$ (such as finite VC dimension, or boundary fragments) or
assumptions on the regularity of the regression function of
classification (plug-in type assumptions). In \cite{nedelec}, minimax
fast rates are stated for finite VC classes of candidates whereas
plug-in type assumptions have been studied in the binary classification
model in \cite{AT} (see also \cite{jaune,yang}). More generally,
\cite{kolt} proposes to consider $\rho>0$ as a complexity parameter
in local Rademacher complexities. It gives general upper bounds
generalizing \eqref{fastrates} and the results of \cite{mammen} and
\cite{AT}. In the present work, a plug-in type complexity assumption
will be considered.

In all these results, empirical risk minimizers appear as good
candidates to reach these fast rates of convergence. Indeed, given a
class of candidates $\mathcal{G}$, a natural way to estimate
$G_K^\star$ is to consider an Empirical Risk Minimization (ERM)
approach. In standard discriminant analysis (e.g., in the free-noise
case considered in \cite{mammen}), the risk $R_K(G)$ in \eqref{bayes}
can be estimated by:\vspace*{-1pt}
%
\begin{equation}
\label{er} R_{n,m}(G)=\frac{1}{2n}\sum
_{j=1}^n \mathbf{1}_{\{ X_j^{(1)}\in
K/G \}} +
\frac{1}{2m}\sum_{j=1}^m
\mathbf{1}_{\{
X_j^{(2)}\in G \}}.
\end{equation}
It leads to an empirical risk minimizer $\hat G_{n,m}$, if it exists,
defined as:\vspace*{-1pt}
%
\begin{equation}
\label{erm} \hat{G}_{n,m}=\arg\min_{G\in\GG}R_{n,m}(G).
\end{equation}

Unfortunately, in the errors-in-variables model, since we observe noisy
samples $Z=X+\epsilon$, the probability densities of the observed
variables w.r.t. the Lebesgue measure are respectively convolution
$(f\cdot\mu)*\eta$ and $(g\cdot\mu)*\eta$, where, for instance, $f\cdot\mu
(x)=f(x)\times\mu(x)$ for all $x\in\mathbb{R}^d$. As a result,
classical ERM principle fails since:\vspace*{-1pt}
\begin{eqnarray*}
&&\frac{1}{2n}\sum_{i=1}^n
\mathbf{1}_{\{ Z_i^{(1)}\in K/G
\}} +\frac{1}{2m}\sum_{i=1}^m
\mathbf{1}_{\{ Z_i^{(2)}\in
G \}} 
\\
&&\quad\mathop{\stackrel{\mathrm{a.s.}} {-\hspace*{-4pt}-\hspace*{-4pt}-
\hspace*{-4pt}\longrightarrow}}_{n,m\rightarrow\infty} \frac
{1}{2} \biggl[\int
_{K/G} (f\cdot\mu)*\eta(x)\,\mathrm{d}x +\int
_G (g\cdot\mu)*\eta (x)\,\mathrm{d}x \biggr]\neq
R_K(G).
\end{eqnarray*}
As a consequence, we add a deconvolution step in the classical ERM
procedure and study the solution of the minimization:\vspace*{-1pt}
\[
\min_{G\in\mathcal{G}}R_{n,m}^\lambda(G),
\]
where $R_{n,m}^\lambda(G)$ is an asymptotically unbiased estimator of
$R_K(G)$. This empirical risk uses kernel deconvolution estimators with
smoothing parameter $\lambda$. It is called deconvolution empirical
risk and will be of the form:
%
\begin{equation}
\label{derm} R_{n,m}^\lambda(G)=\frac{1}{2n} \sum
_{j=1}^n h_{K/G,\lambda
}
\bigl(Z_j^{(1)}\bigr)+\frac{1}{2m}\sum
_{j=1}^m h_{G,\lambda}\bigl(Z_j^{(2)}
\bigr),
\end{equation}
where the $h_{G,\lambda}(\cdot)$ are deconvoluted versions of
indicator functions used in classical ERM for direct observations (see
Section~\ref{section:upperbounds} for details).

In this contribution, we would like to describe as precisely as
possible the influence of the error $\epsilon$ on the classification
rates\vadjust{\goodbreak} of convergence and the presence of fast rates. Our aim is to use
the asymptotic theory of empirical processes in the spirit of \cite
{vdg} (see also \cite{wvdv}) when dealing with the deconvolution
empirical risk \eqref{derm}. To this end, we study in details the
complexity of the class of functions $\{h_{G,\lambda},G\in\mathcal
{G}\}$, given the explicit form of functions $h_{G,\lambda}$. This
complexity is related to the imposed complexity over $\mathcal{G}$.

We establish lower and upper bounds and discuss the performances of
this deconvolution ERM estimator under a plug-in complexity assumption.
As mentioned earlier, different complexity assumptions have been
developed in the last decades. The boundary fragment regularity,
considered by, e.g., \cite{korostelevtsybakov,mammen} is the core of a
future work.

We point out that the definition of the empirical risk \eqref{derm}
leads to a new and interesting theory of risk bounds detailed in
Section~\ref{section:upperbounds} for discriminant analysis. In particular, parameter $\lambda
$ has to be calibrated to reach a bias/variance trade-off in the
decomposition of the excess risk. Related ideas have been recently
introduced in \cite{klemela} in the Gaussian white noise model and
density estimation setting for more general linear inverse problems
using singular values decomposition. In our framework, up to our
knowledge, the only minimax result is \cite{mme} which gives minimax
rates in Hausdorff distance for manifold estimation in the presence of
noisy variables. \cite{delgib} gives also consistency and limiting
distribution for estimators of boundaries in deconvolution problems,
but no minimax results are proposed. In the free-error case, we can
also apply this methodology. In this case, the empirical risk is given
by the estimation of $f$ and $g$ using simple kernel density
estimators. This idea has been already mentioned in \cite{vapnik2000}
in the general learning context and called Vicinal Risk Minimization
(see also \cite{vrm}). However, even in pattern recognition and in the
direct case, up to our knowledge, there is no asymptotic rates of
convergence for this empirical minimization principle.

In this contribution, a classifier $G$ is always identified with a
subset of $\R^d$. Our aim is to mimic the set $G_K^\star$ from the
noisy observations (\ref{eq:modelb}). In particular, we aim at
understanding the relationship between the spatial position of an input
$X\in\R^d$ and its affiliation to one of the candidate densities. For
this purpose, we give a deconvolution strategy to minimize the excess
risk \eqref{risk}. This problematic falls into the general problem of
prediction with measurement errors (see \cite{Carroll_Delaigle_Hall}).
This is the classification counterpart of the more extensively studied
model of regression with errors-in-variables (see \cite{fan2} or more
recently \cite{meister}). It is important to note that one could
alternatively try to provide the best classifier for a noisy input $Z$.
In this case, we are faced to a direct problem which is in some sense
already treated in \cite{mammen}. However, it could be interesting to
compare the performances of the two different approaches.

At this step, remark that similar problems have been considered in the
test theory. Indeed, if we deal with a new incoming (noise free)
observation $X$ having density $f_X$, our aim is exactly to test one of
the following `inverse' hypotheses:
%
%
\begin{equation}
\label{eq:pb_inverse} H_0^{IP}\dvtx f_X=f, \quad
\mbox{against}\quad H_1^{IP}\dvtx f_X=g.
\end{equation}
However, we do not set any kind of order (null and alternative) between
$H_0$ and $H_1$. The risk $R_K(G)$ is then related to the sum of the
first and second kind error. Alternatively, if we deal with a noisy
input $Z$ having density $(f_X\cdot\mu)*\eta$, this would correspond to test:
%
%
\begin{equation}
\label{eq:pb_direct} H_0^{DP}\dvtx (f_X\cdot\mu)*
\eta=(f\cdot\mu)*\eta, \quad\mbox{against}\quad H_1^{DP}
\dvtx (f_X\cdot\mu)*\eta=(g\cdot\mu)*\eta.
\end{equation}
A natural question then arises: are the both problems
(\ref{eq:pb_inverse}) and (\ref{eq:pb_direct}) equivalent or comparable?
This question has already been addressed in \cite{test2} or \cite
{test3} in a slightly different setting. This could be the core of a
future work, but it requires the preliminary study provided in these papers.

Finally, for practical motivation, we can refer to the monograph of
Meister \cite{meister} for particular models with measurement errors, such as
in medicine, econometry or astronomy. In the specific context of
classification, we met two explicit examples. The first one is an
example in oncology where we try to classify the evolution of cancer
thanks to medical images (like MRI or X-ray). These images are noisy
due to the data collection process or the interpretation of the
practitioner. The second example comes from meteorology where the
weather forecaster wants to predict the future raining day thanks to
measures such as rain gauge or barometer (which have well-studied
random errors).

The paper is organized as follows. In Section~\ref{section:lowerbounds},
the model assumptions are explicited and an
associated lower bound is stated. This lower bound generalizes to the
indirect case the well-known lower bound of \cite{AT} established in
classification. Deconvolution ERM attaining these rates are presented
in Section~\ref{section:upperbounds}. We also consider in this section
standard kernel estimators, which allow to construct a new minimax
optimal procedure in the direct case. A brief discussion and some
perspectives are gathered in Section~\ref{section:conclusion} while
Section~\ref{section:proofs} is dedicated to the proofs of the main results.

\section{Lower bound}
\label{section:lowerbounds}

\subsection{Model setting}\label{sec2.1}

In this section, we detail some common assumptions (complexity and
margin) on the pair $(f,g)$. We then propose a lower bound on the
corresponding minimax rates.

First of all, given a set $G \subset K$, simple algebra indicates that
the excess risk $R_K(G)-R_K(G_K^\star)$ can be written as:
\[
R_K(G)-R_K\bigl(G_K^\star
\bigr)=\tfrac{1}{2}d_{f,g}\bigl(G,G_K^\star
\bigr),
\]
where the pseudo-distance $d_{f,g}$ over subsets of $K \subset\R^d$
is defined as:
\[
d_{f,g}(G_1,G_2)=\int_{G_1\Delta G_2}
\llvert f-g \rrvert \,\mathrm{d}Q,
\]
and $G_1\Delta G_2=[G_1^c\cap G_2]\cup[G_2^c\cap G_1]$ is the
symmetric difference between two sets $G_1$ and~$G_2$. In this context,
there is another natural way of measuring the accuracy of a decision
rule $G$ through the quantity:
\[
d_{\Delta}\bigl(G,G_K^\star\bigr)=\int
_{G\Delta G_K^\star}\,\mathrm{d}Q,
\]
where $d_{\Delta}$ defines also a pseudo-distance on the subsets of $K
\subset\R^d$.

In this paper, we are interested in the minimax rates associated to
these pseudo-distances. In other words, given a class $\mathcal{F}$,
one would like to quantify as precisely as possible the corresponding
minimax risks defined as
\[
\inf_{\hat G_{n,m}} \sup_{(f,g)\in\mathcal{F}}
\E_{f,g}d_{\square
}\bigl(\hat G_{n,m},G_K^\star
\bigr),
\]
where the infimum is taken over all possible estimators of $G_K^\star$
and $d_{\square}$ stands for $d_{f,g}$ or $d_\Delta$ following the
context. In particular, we will exhibit classification rules $\hat
G_{n,m}$ attaining these rates. In order to obtain a satisfying study
of the minimax rates mentioned above, one needs to detail the
considered classes $\mathcal{F}$. Such a class expresses some
conditions on the pair $(f,g)$. They are often separated into two
categories: margin and complexity assumptions.

A first condition is the well-known Margin assumption. It has been
introduced in discriminant analysis (see \cite{mammen}) as follows.

\begin{MarginAssumption*}
There exists positive constants
$t_0,c_2,\alpha\geq0$ such that for $0<t<t_0$:
%
\begin{equation}
\label{ma} Q\bigl\{x\in K\dvtx \bigl\llvert f(x)-g(x) \bigr\rrvert \leq t
\bigr\}\leq c_2t^{\alpha}.
\end{equation}
\end{MarginAssumption*}

This assumption is related to the behavior of $\llvert
f-g
\rrvert $ at the boundary
of $G_K^\star$. It may give a variety of minimax fast rates of
convergence which depends on the margin parameter $\alpha$. A~large
margin corresponds to configurations where the slope of $\llvert
f-g
\rrvert $ is high
at the boundary of $G_K^\star$. The most favorable case arises when
the margin $\alpha=+\infty$. In such a situation, $f-g$ has a
discontinuity at the boundary of $G_K^\star$.

From a practical point of view, this assumption provides a precise
description of the interaction between the pseudo distance $d_{f,g}$
and $d_{\Delta}$. In particular, it allows a control of the variance
of the empirical processes involved in the upper bounds, thanks to
Lemma~2 in \cite{mammen}. More general assumptions of this type can be
formulated (see, for instance, \cite{empimini} or \cite{kolt}) in a
more general statistical learning context.

For the sake of convenience, we will require in the following an
additional assumption on the noise $\epsilon$. We assume in the sequel
that $\epsilon= (\epsilon_1,\ldots,\epsilon_d)'$ admits a bounded
density $\eta$ with respect to the Lebesgue measure satisfying:
%
%
\begin{equation}
\label{eq:bruitindep} \eta(x) = \prod_{i=1}^d
\eta_i(x_i)\qquad \forall x\in\R^d.
\end{equation}
In other words, the entries of the vector $\epsilon$ are independent.
The assumption below describes the difficulty of the considered
problems. It is often called the ordinary smooth case in the inverse
problem literature.

\begin{NoiseAssumption*}
There exist $(\beta_1,\ldots
,\beta_d)'\in\R_+^d$ and $\mathcal{C}_1,\mathcal{C}_2,\mathcal
{C}_3$ positive constants such that for all $i\in\{1,\ldots, d
\}$, $\beta_i>1/2$,
\[
\mathcal{C}_1 \llvert t \rrvert ^{-\beta_i} \leq\bigl\llvert
\mathcal{F}[\eta_i](t) \bigr\rrvert \leq\mathcal{C}_2
|t|^{-\beta_i}, \quad\mbox{and}\quad \biggl\llvert \frac{\mathrm{d}}{\mathrm{d}t}
\mathcal{F}[\eta_i](t) \biggr\rrvert \leq\mathcal{C}_3
|t|^{-\beta_i} \qquad\mbox{as } |t|\to+\infty,
\]
where $\mathcal{F}[\eta_i]$ denotes the Fourier transform of $\eta
_i$. Moreover, we assume that $\mathcal{F}[\eta_i](t) \neq0$ for
all $t\in\R$ and $i\in\{1,\ldots, d \}$.
\end{NoiseAssumption*}

Classical results in deconvolution (see, e.g., \cite{fan,fan2} or
\cite{Butucea} among others) are stated for $d=1$. Two
different settings are then distinguished concerning the difficulty of
the problem which is expressed through the shape of $\mathcal{F}[\eta
]$. One can consider alternatively the case where $\mathcal{C}_1
|t|^{-\beta} \leq|\mathcal{F}[\eta](t)| \leq\mathcal{C}_2
|t|^{-\beta}$ as $|t|\rightarrow+\infty$, which yet corresponds to
mildly ill-posed inverse problem or $\mathcal{C}_1 \mathrm{e}^{-\gamma
|t|^\beta} \leq|\mathcal{F}[\eta](t)| \leq\mathcal{C}_2
\mathrm{e}^{-\gamma|t|^\beta}$, $\gamma>0$ as $|t|\rightarrow+\infty$ which
leads to a severely ill-posed inverse problem. This last setting
corresponds to a particularly difficult problem and is often associated
to low minimax rates of convergence.

In this contribution, we only deal with $d$-dimensional mildly
ill-posed deconvolution problems. For the sake of brevity, we do not
consider severely ill-posed inverse problems or possible intermediates
(e.g., a combination of polynomial and exponential decreasing
functions). Nevertheless, the rates in these cases could be obtained
through the same steps.

The Margin assumption is `structural' in the sense that it describes
the difficulty to distinguish an observation having density $f$ from an
other with density $g$. In order to provide a complete study, one also
needs to set an assumption on the difficulty to find $G_K^\star$ in a
possible set of candidates, namely a complexity assumption. In the
classification framework, two different kinds of complexity assumptions
are often introduced in the literature. The first kind concerns the
regularity of the boundary of the Bayes classifier. Indeed, our aim is
to estimate $G_K^\star$, which yet corresponds to a nonparametric set
estimation problem. In this context, it seems natural to traduce the
difficulty of the learning process by condition on the shape of
$G_K^\star$. Another way to describe the complexity of the problem is
to impose condition on the regularity of the underlying densities $f$
and $g$. Such kind of condition is originally related to plug-in
approaches and will be the investigated framework. Remark that these
two assumptions are quite different and are convenient for distinct
problems. In particular, a set $G_K^\star$ with a smooth boundary is
not necessarily associated to smooth densities, and vice-versa.

In the rest of this section, lower bounds for the associated minimax
rates of convergence are stated in the noisy setting. Corresponding
upper bounds are presented and discussed in Section~\ref
{section:upperbounds}.

\subsection{Lower bound for the Plug-in assumption}

The Plug-in assumption considered in this paper is related to the
regularity of the function $f-g$, expressed in terms of H\"older
spaces. It corresponds to the same kind of assumption as in \cite{AT}
for classification.

Given $\gamma,L>0$, $\Sigma(\gamma,L)$ is the class of isotropic H\"
older continuous functions $\nu$ having continuous partial derivatives
up to order $\lfloor\gamma\rfloor$, the maximal integer strictly
less than $\gamma$ and such that:
\[
\bigl |\nu(y)-p_{\nu,x}(y)\bigr |\leq L\|x-y\|^\gamma,\qquad\forall x,y\in
\R^d,
\]
where $p_{\nu,x}$ is the Taylor polynomial of $\nu$ at order $\lfloor
\gamma\rfloor$ at point $x$ and $\|\cdot\|$ stands for the Euclidean norm
on $\mathbb{R}^{d}$.

\begin{PluginAssumption*}
There exist positive constants
$\gamma$ and $L$ such that $f-g \in\Sigma(\gamma,L)$.
\end{PluginAssumption*}

We then call $\mathcal{F}_{\mathrm{plug}}(Q)$ the set of all pairs
$(f,g)$ satisfying both the  {Margin} (with respect to $Q$) and
the  {Plug-in} assumptions, since the previous assumption is
often associated to plug-in rules in the statistical learning
literature. The following theorem proposes a lower bound for the noisy
discriminant analysis problem in such a setting.

%
\begin{them}
\label{thm:lbplugin}
Suppose that the Noise assumption is satisfied. Then, there exists a
measure $Q_0$ such that for all $\alpha\leq1$,
\[
\liminf_{n,m \to+\infty} \inf_{\hat G_{n,m}} \sup
_{(f,g)\in
\mathcal{F}_{\mathrm{plug}}(Q_0)} (n\wedge m)^{\tau_d(\alpha,\beta
,\gamma)} \E_{f,g}
d_{\square}\bigl(\hat G_{n,m},G_K^{\star}
\bigr) > 0,
\]
where the infinimum is taken over all possible estimators of the set
$G_K^{\star}$ and
\[
\tau_d(\alpha,\beta,\gamma) =\cases{ %
\displaystyle\frac{ \gamma\alpha}{ \gamma(2+\alpha) +d +2\sum_{i=1}^{d}\beta
_i } &
\quad$\mbox{for }d_{\square}=d_\Delta$,\vspace*{2pt}
\cr
\\
\displaystyle\frac{ \gamma(\alpha+1)}{ \gamma(2+\alpha)+d+ 2\sum_{i=1}^{d}\beta_i }&\quad $\mbox{for } d_{\square}=d_{f,g}$. }
\]
\end{them}

Remark that we obtain exactly the same lower bounds as \cite{AT} in
the direct case, which yet corresponds to the situation where $\beta
_j=0$ for all $j\in\{1, \ldots, d \}$.

In the presence of noise in variables, the rates obtained in Theorem~\ref{thm:lbplugin} are slower. The price to pay is an additional term
of the form:
\[
2\sum_{i=1}^{d}\beta_i.
\]
This term clearly connects the difficulty of the problem to the tail
behavior of the characteristic function of the noise distribution. This
price to pay is already known in density estimation, regression with
errors in variables or goodness-of-fit testing. Last step is to get a
corresponding upper bound to validate this lower bound in the presence
of noise in variables.

Remark that this lower bound is valid only for $\alpha\leq1$. This
restriction appears for some technical reasons in the proof (see
Section~\ref{section:proofs}). The main difficulty here is to use
standard arguments from lower bounds in classification (see \cite
{audibert2004,AT}) in this deconvolution setting. More precisely, we
have to take advantage of the Noise assumption, related to the Fourier
transform of the noise distribution $\eta$. To this end, we use in the
proof of Theorem~\ref{thm:lbplugin} an algebra based on standard
Fourier analysis tools, and we have to consider sufficiently smooth
objects. As a consequence in the lower bounds, we can check the Margin
assumption only for values of $\alpha\leq1$. Nevertheless, we
conjecture that this restriction is only due to technical reasons and
that our result remains pertinent for all $\alpha\geq0$. In
particular, an interesting direction is to consider a wavelet basis
which provides an isometric wavelet transform in $L^2$ in order to
obtain the desired lower bound in the general case.

The measure $Q_0$ that we mention in Theorem~\ref{thm:lbplugin} is
explicitly constructed in the proof. For the sake of convenience, the
construction of this measure is not reproduced here (we refer to
Section~\ref{s:proof_lb} for an interested reader).


\section{Upper bounds}
\label{section:upperbounds}

\subsection{Estimation of \texorpdfstring{$G_K^\star$}{G*K}}

In the free-noise case ($\epsilon_j^{(i)}=(0,\ldots, 0)$ for all
$j\in\{1,\ldots,n \}$, $i\in\{1,2 \}$), we
deal with two samples $(X_1^{(1)},\ldots,X_n^{(1)})$,
$(X_1^{(2)},\ldots,X_m^{(2)})$ having respective densities $f$ and
$g$. A~standard way to estimate $G_K^\star=\{ x\in K\dvtx  f(x) \geq
g(x) \}$ is to estimate $R_K(\cdot)$ thanks to the data. For all
$G\subset K$, the risk $R_K(G)$ can be estimated by the empirical risk
defined in \eqref{er}. Then the Bayes classifier $G_K^{\star}$ is
estimated by $\hat G_{n,m}$ defined as a minimizer of the empirical
risk \eqref{er} over a given family of sets $\GG$. We know for
instance from \cite{mammen} that the estimator $\hat{G}_{n,m}$
reaches the minimax rates of convergence in the direct case when
$\mathcal{G}=\GG(\gamma,L)$ corresponds to the set of boundary
fragments with $\gamma>d-1$. For larger set $\GG(\gamma,L)$, as
proposed in \cite{mammen}, the minimization can be restricted to a
$\delta$-net of $\GG(\gamma,L)$. With an additional assumption over
the approximation power of this $\delta$-net, the same minimax rates
can be achieved in a subset of $\GG(\gamma,L)$.

If we consider complexity assumptions related to the smoothness of
$f-g$, we can show easily with \cite{AT} that an hybrid plug-in/ERM
estimator reaches the minimax rates of convergence of \cite{AT} in the
free-noise case. The principle of the method is to consider the
empirical minimization \eqref{erm} over a particular class $\GG$
based on plug-in type decision sets. More precisely, following \cite
{AT} for classification, we can minimize in the direct case the
empirical risk over a class $\GG$ of the form:
\[
\GG=\bigl\{\{f-g\geq0\},f-g\in\mathcal{N}_{n,m}\bigr\},
\]
where $\mathcal{N}_{n,m}$ is a well-chosen $\delta$-net. With such a
procedure, minimax rates can be obtained with no restriction over the
parameter $\gamma$, $\alpha$ and $d$.

In noisy discriminant analysis, ERM estimator \eqref{erm} is no longer
available as mentioned earlier.
Hence, we have to add a deconvolution step to the classical ERM
estimator. In this context, we can construct a deconvolution kernel,
provided that the noise has a nonnull Fourier transform, as expressed
in the {Noise assumption}. Such an assumption is rather
classical in the inverse problem literature (see, e.g., \cite{fan,Butucea}
or \cite{meister}).

Let $\mathcal{K}=\prod_{j=1}^d \mathcal{K}_j\dvtx \R^d \to\R$ be a
$d$-dimensional function defined as the product of $d$ unidimensional
functions $\mathcal{K}_j$. The properties of $\mathcal{K}$ leading to
satisfying upper bounds will be made precise later on. Then, if we denote by
$\lambda=(\lambda_1,\ldots,\lambda_d)$ a set of (positive)
bandwidths and by $\FF[\cdot]$ the Fourier transform, we define
$\mathcal{K}_\eta$ as:
%
%
\begin{eqnarray}
\label{dk} \mathcal{K}_{\eta} \dvtx \R^d &\to&\R,
\nonumber
\\[-8pt]
\\[-8pt]
t &\mapsto&\mathcal{K}_\eta(t) = \FF^{-1} \biggl[
\frac{\FF
[\mathcal{K}](\cdot)}{\FF[\eta](\cdot/\lambda)} \biggr](t).
\nonumber
\end{eqnarray}
In this context, for all $G\subset K$, the risk $R_K(G)$ can be
estimated by
\[
R^{\lambda}_{n,m}(G)=\frac{1}{2} \Biggl[\frac{1}{n}
\sum_{j=1}^nh_{K/ G,\lambda}
\bigl(Z_j^{(1)}\bigr)+\frac{1}{m}\sum
_{j=1}^mh_{G,\lambda}\bigl(Z_j^{(2)}
\bigr) \Biggr], 
\]
where for a given $z\in\R^d$:
%
\begin{equation}
\label{hG} h_{G,\lambda}(z)= \int_{G}
\frac{1}{\lambda} \mathcal{K}_\eta \biggl(\frac{z-x}{\lambda} \biggr)\,
\mathrm{d}x.
\end{equation}
In the following, we study ERM estimators defined as:
%
\begin{equation}
\label{derm2} \hat{G}_{n,m}^\lambda=\arg\min
_{G\in\GG}R_{n,m}^\lambda(G),
\end{equation}
where parameter $\lambda=(\lambda_1,\ldots,\lambda_d)\in\R_+^{d}$
has to be chosen explicitly. Functions $h_{G,\lambda}$ in equation
\eqref{hG} are at the core of the upper bounds. In particular,
following the pioneering's works of Vapnik (see \cite{vapnik2000}), we
have for $R^\lambda_K(\cdot):=\E R_{n,m}^{\lambda}(\cdot)$:
%
\begin{eqnarray}
\label{noisyvapnik} R_K\bigl(\hat{G}^\lambda_{n,m}
\bigr)-R_K\bigl(G_K^\star\bigr) &\leq&
R_K\bigl(\hat {G}^\lambda_{n,m}
\bigr)-R_{n,m}^\lambda\bigl(\hat{G}^\lambda
_{n,m}\bigr)+R_{n,m}^\lambda\bigl(G_K^\star
\bigr)-R_K\bigl(G_K^\star\bigr)
\nonumber
\\
&\leq& R^\lambda_K\bigl(\hat{G}^\lambda_{n,m}
\bigr)-R_{n,m}^\lambda\bigl(\hat {G}^\lambda_{n,m}
\bigr)+R_{n,m}^\lambda\bigl(G_K^\star
\bigr)-R^\lambda_K\bigl(G_K^\star
\bigr)
\nonumber
\\[-8pt]
\\[-8pt]
&&{} +\bigl(R_K-R_K^\lambda\bigr) \bigl(
\hat{G}^\lambda_{n,m}\bigr)-\bigl(R_K-R^\lambda
_K\bigr) \bigl(G_K^\star\bigr)
\nonumber
\\
&\leq& \sup_{G\in\GG}\bigl |R_K^\lambda-R_{n,m}^\lambda\bigr |
\bigl(G,G_K^\star\bigr) +\sup_{G\in\GG}\bigl |R_K^\lambda-R_K\bigr |
\bigl(G,G_K^\star\bigr),
\nonumber
\end{eqnarray}
where we write for concision for any $G,G'\subset K$:
\[
\bigl |R_K^\lambda-R_{n,m}^\lambda\bigr |
\bigl(G,G'\bigr)=\bigl |R_K^\lambda(G)-R_K^\lambda
\bigl(G'\bigr)-R_{n,m}^\lambda(G)+R_{n,m}^\lambda
\bigl(G'\bigr)\bigr |,
\]
and similarly:
\[
\bigl |R_K^\lambda-R_{K}\bigr |\bigl(G,G'
\bigr)=\bigl |R_K^\lambda(G)-R_K^\lambda
\bigl(G'\bigr)-R_{K}(G)+R_{K}
\bigl(G'\bigr)\bigr |.
\]
As a result, to get risk bounds, we have to deal with two opposing
terms, namely a so-called variability term:
%
\begin{equation}
\label{variance} \sup_{G\in\GG}\bigl |R_K^\lambda-R_{n,m}^\lambda\bigr |
\bigl(G-G_K^\star\bigr),
\end{equation}
and a bias term (since $\E R_{n,m}^\lambda(G)\neq R_K(G)$) of the form:
%
\begin{equation}
\label{bias} \sup_{G\in\GG}\bigl |R_K^\lambda-R_K\bigr |
\bigl(G-G_K^\star\bigr).
\end{equation}
The variability term \eqref{variance} gives rise to the study of
increments of empirical processes. In this work, this control is based
on entropy conditions and uniform concentration inequalities. It is
inspired by results presented for instance in \cite{wvdv} or \cite
{vdg}. The main novelty here is that in the noisy case, empirical
processes are indexed by a class of functions which depends on the
smoothing parameter $\lambda$. The bias term \eqref{bias} is
controlled by taking advantages of the properties of $\mathcal{G}$ and
of the assumptions on the kernel $\mathcal{K}$. Indeed, it can be
related to the standard bias term in nonparametric density estimation
and can be controlled using smoothness assumptions of plug-in type.
This bias term is inherent to the estimation procedure and its control
is a cornerstone of the upper bounds.

The choice of $\lambda$ will be a trade-off between the two opposing
terms \eqref{variance} and \eqref{bias}. Small $\lambda$ leads to
complex functions $h_{G,\lambda}$ and blasts the variance term whereas
\eqref{bias} vanishes when $\lambda$ tends to zero. The kernel
$\mathcal{K}$ has to be chosen in order to take advantage of the
different conditions on $G_K^\star$. This choice will be operated
according to the following definition.

\begin{Definition*}
We say that $\mathcal{K}$ is a kernel of order
$l\in\mathbb{N}^*$ if and only if:
\begin{itemize}
\item$\int\mathcal{K}(u)\,\mathrm{d}u=1$.
\item$\int u_j^k\mathcal{K}(u)\,\mathrm{d}u=0$ $\forall k=1,\ldots, l$,
$\forall j=1,\ldots, d$.
\item$\int|u_j|^{l+1}|\mathcal{K}(u)|\,\mathrm{d}u<\infty$, $\forall j=1,\ldots, d$.
\end{itemize}
\end{Definition*}

In addition to this definition, we will require the following
assumption on the kernel $\mathcal{K}$ which appears in \eqref{dk}.


\begin{KernelAssumption*}
The kernel $\mathcal{K}$ is such
that $\mathcal{F}[\mathcal{K}]$ is bounded and compactly supported.
\end{KernelAssumption*}

The construction of kernels of order $l$ satisfying the Kernel
assumption could be managed using for instance the so-called Meyer
wavelet (see \cite{mallat}).

The following subsection intent to study deconvolution ERM estimator
\eqref{derm2} and gives asymptotic fast rates of convergence. It
validates the lower bounds of Theorem~\ref{thm:lbplugin}.

%
%
%

\subsection{Upper bound}
\label{section:upper_bound}

For all $\delta>0$, using the notion of entropy (see, for instance,
\cite{wvdv}) for H\"olderian function on compact sets, we can find a
$\delta$-network $\mathcal{N}_{\delta}$ on $\Sigma(\gamma,L)$ such that:
\begin{itemize}
\item$\log(\operatorname{card}(\mathcal{N}_\delta)) \leq A \delta
^{-d/\gamma}$,
\item For all $h_0\in\Sigma(\gamma,L)$, we can find $h\in\mathcal
{N}_{\delta}$ such that $\| h-h_0 \|_\infty\leq\delta$.
\end{itemize}
In the following, we associate to each $\nu:=f-g \in\Sigma(\gamma
,L)$, a set $G_\nu= \{ x\in K \dvtx \nu(x)\geq0 \}$ and define
the ERM estimator as:
%
%
\begin{equation}
\label{eq:sieve} \hat G_{n,m} = \operatorname{arg} \min
_{\nu\in\mathcal{N}_\delta} R_{n,m}^\lambda(G_\nu),
\end{equation}
where $\delta=\delta_{n,m}$ has to be chosen carefully. This
procedure has been introduced in the direct case by \cite{AT} and
referred to as an hybrid plug-in/ERM procedure. The following theorem
describes the performances of $\hat G_{n,m}$.

%
\begin{them}
\label{thm:upplugin}
Let $\hat G_{n,m}$ the set introduced in (\ref{eq:sieve}) with
\begin{eqnarray*}
\lambda_j &=& (n\wedge m)^{-\trup{1}{(\gamma(2+\alpha) +2 \sum
_{i=1}^d \beta_i + d)}},\qquad \forall j\in\{1,
\ldots,d \}, \quad\mbox{and}
\\
 \delta&=& \delta_{n,m} = \biggl(
\frac{\prod_{i=1}^d
\lambda_i^{-\beta_i}}{\sqrt{n\wedge m}} \biggr)^{\trup{2}{(d/\gamma
+2+\alpha)}}.
\end{eqnarray*}
Given some $\sigma$-finite measure $Q$, suppose $(f,g)\in\mathcal
{F}_{\mathrm{plug}}(Q)$ and the Noise assumption is satisfied with
$\beta_i>1/2$, $\forall i=1,\ldots, d$. Consider a kernel $\mathcal
{K}_\eta$ defined as in \eqref{dk} where $\mathcal{K}=\prod_{j=1}^d\mathcal{K}_j$ is a kernel of order $\lfloor\gamma\rfloor$,
which satisfies the Kernel assumption. Then, for all real $\alpha\geq
0$, if $Q$ is the Lebesgue measure:
\[
\lim_{n,m \to+\infty} \sup_{(f,g)\in\mathcal{F}_{\mathrm
{plug}}(Q)} (n\wedge
m)^{\tau_d(\alpha,\beta,\gamma)} \E_{f,g} d_\square\bigl(\hat
G_{n,m},G_K^{\star}\bigr) < +\infty,
\]
where
\[
\tau_d(\alpha,\beta,\gamma) =\cases{ %
\displaystyle\frac{ \gamma\alpha}{ \gamma(2+\alpha) +d +2\sum_{i=1}^{d}\beta
_i } &
\quad$\mbox{for }d_\square=d_\Delta$,\vspace*{2pt}
\cr
\displaystyle\frac{ \gamma(\alpha+1)}{ \gamma(2+\alpha)+d+ 2\sum_{i=1}^{d}\beta_i } &
\quad$\mbox{for } d_\square=d_{f,g}$. }
\]
Moreover, if $Q(x)=\mu(x)\,\mathrm{d}x$, the same upper bounds hold provided that
$\mu\in\Sigma(\gamma,L)$ and that
$\min_{x\in K} \mu(x) \geq c_0$ for some $c_0>0$.
\end{them}

Theorem~\ref{thm:upplugin} validates the lower bounds of Theorem~\ref
{thm:lbplugin}. Deconvolution ERM are minimax optimal over the class
$\mathcal{F}_{\mathrm{plug}}$. These optimal rates are characterized
by the tail behavior of the characteristic function of the error
distribution $\eta$. We only consider the ordinary smooth case whereas
straightforward modifications lead to slow rates of convergence in the
super-smooth case.

Here, fast rates (i.e., faster than $1/\sqrt{n}$) are pointed out when
$\alpha\gamma>d+2\sum\beta_i$. This result is comparable to \cite
{AT}, where fast rates are proposed when $\alpha\gamma>d$. However,
it is important to stress that large values of both $\alpha$ and
$\gamma$ correspond to restrictive situations. In this case, the
margin parameter is high whereas the behavior of $f-g$ is smooth, which
seems to be contradictory (see the related discussion in \cite{AT}).

If $Q$ is not the Lebesgue measure, $\mu$ has to be lower bounded by a
constant $c_0>0$ which appears in the upper bound. This assumption can
be relaxed to recover the case of Theorem~\ref{thm:lbplugin}.

For the sake of concision, we do not study plug-in rules in this paper.
Such algorithms are characterized by classifiers of the form
\[
\tilde G_{n,m} = \bigl\{ x\in K, \tilde f_n(x) -\tilde
g_m(x) \geq0 \bigr\},
\]
where $\tilde f_n -\tilde g_m$ is an (optimal) estimator of the
function $f-g$. The performances of such kind of methods have been
investigated by \cite{AT} in the binary classification model. We also
mention for instance \cite{Goldstein_Messer1992} or \cite
{Bickel_Ritov2003} for contributions in a more general framework.

Nevertheless, we point out that the choice of $\lambda$ in Theorem~\ref{thm:upplugin} is the trade-off between the variability term
\eqref{variance} and the bias term \eqref{bias}. It is important to
note that this asymptotic for $\lambda$ is not the optimal choice in
the problem of deconvolution estimation of $f-g\in\Sigma(\gamma,L)$
thanks to noisy data. Here the bandwidth depends on the margin
parameter $\alpha$ and optimizes the classification excess risk bound.
It highlights that the estimation procedure \eqref{eq:sieve} is not a
plug-in rule but an hybrid ERM/Plug-in estimator as in \cite{AT}.

Finally, this deconvolution ERM appears to be minimax optimal when we
deal with noisy data such that $\beta_i>\frac{1}{2}$, $\forall
i=1,\ldots, d$. A natural question is to extend these results to the
direct case where $\beta_i=0$, $\forall i=1,\ldots, d$. Moreover, the
minimax optimality of this procedure depends on the choice of $\lambda
$ in Theorem~\ref{thm:upplugin}. In the following subsection, we deal
with a similar approach in the direct case, using standard kernel
estimators instead of deconvolution kernel estimators. Interestingly in
this situation, the choice of $\lambda$ is not crucial to derive
optimal rates of convergence.

\subsection{Upper bound in the free-noise case}

In the free-noise setting, direct observations $X_j^{(1)}$, $j=1,\ldots
,n$ and $X^{(2)}_j$, $j=1,\ldots, m$ are available. In this case, we
can construct an estimation procedure based on \eqref{eq:sieve} where
a standard kernel estimator is used instead of a deconvolution kernel
estimator. Following the noisy setting, we define in the direct case
$\tilde G^\lambda_{n,m}$ as follows:
%
%
\begin{equation}
\label{sieve2} \tilde G^\lambda_{n,m} = \operatorname{arg}
\min_{\nu\in\mathcal
{N}_\delta} \tilde R_{n,m}^\lambda(G_\nu),
\end{equation}
where here $\tilde R^\lambda_{n,m}(G)$ is an estimator of $R_K(G)$
defined as:
\[
\tilde R^\lambda_{n,m}(G)=\frac{1}{2} \Biggl[
\frac{1}{n}\sum_{j=1}^n\tilde
h_{K/ G,\lambda}\bigl(X_j^{(1)}\bigr)+\frac{1}{m}
\sum_{j=1}^m\tilde h_{G,\lambda}
\bigl(X_j^{(2)}\bigr) \Biggr], 
\]
where for a given kernel $\mathcal{K}$:
\[
\tilde h_{G,\lambda}(z)= \int_{G}
\frac{1}{\lambda} \mathcal{K} \biggl(\frac{z-x}{\lambda} \biggr)\,\mathrm{d}x.
\]
The following theorem describes the performances of $\tilde G^\lambda_{n,m}$.

%
\begin{cor}
\label{cor}
Let $\mathcal{F} = \mathcal{F}_{\mathrm{plug}}(Q)$ and $\tilde
G^\lambda_{n,m}$ the set introduced in (\ref{sieve2}) with
\[
\lambda_j \leq(n\wedge m)^{-\trup{1}{(\gamma(2+\alpha) + d)}}, \qquad\forall j\in\{1,
\ldots,d \}, \quad\mbox{and}\quad \delta=\delta _{n,m} = \biggl(
\frac{1}{\sqrt{n\wedge m}} \biggr)^{\trup
{2}{(d/\gamma+2+\alpha)}}.
\]
Consider a kernel $\mathcal{K}=\prod_{j=1}^d\mathcal{K}_j$ of order
$\lfloor\gamma\rfloor$ satisfying the Kernel assumption. Then, if
$Q$ is the Lebesgue measure, for any real $\alpha\geq0$:
\[
\lim_{n,m \to+\infty} \sup_{(f,g)\in\mathcal{F}_{\mathrm
{plug}}(Q)} (n\wedge
m)^{\tau_d(\alpha,\gamma)} \E d_\square\bigl(\tilde G^\lambda_{n,m},G_K^{\star}
\bigr) < +\infty,
\]
where
\[
\tau_d(\alpha,\gamma) =\cases{ %
\displaystyle\frac{ \gamma\alpha}{ \gamma(2+\alpha) +d } &
\quad$\mbox{for }d_\square =d_\Delta$,\vspace*{2pt}
\cr
\\
\displaystyle\frac{ \gamma(\alpha+1)}{ \gamma(2+\alpha)+d } &\quad$\mbox{for } d_\square=d_{f,g}$. }
\]
Moreover, if $Q(x)=\mu(x)\,\mathrm{d}x$, the same upper bounds holds provided
that $\mu\in\Sigma(\gamma,L)$ and that
$\min_{x\in K} \mu(x) \geq c_0$ for some $c_0>0$.
\end{cor}

These rates correspond to the lower bound of Theorem~\ref
{thm:lbplugin} for $\beta_j=0$, $\forall j=1,\ldots, d$ (see also
\cite{AT}). As a result, \eqref{sieve2} provides a new procedure
which reaches the minimax optimality in classification. Some remarks
are in order.

The choice of $\lambda$ in Corollary~\ref{cor} is not standard. It seems that
if $\lambda$ is small enough, the ERM procedure \eqref{sieve2} is
minimax. This result can be explain as follows. Here, $\lambda$ is not
a trade-off between two opposing terms. In the control of the
variability term, it appears that with a good choice of $\mathcal{K}$,
the variability term does not depend on the bandwidth $\lambda$ of the
kernel. As a result, we only need to control the bias term with a small
bandwidth.

This property can also be interpreted heuristically as follows. It is
clear that the estimation procedure \eqref{sieve2} with kernel
estimator $\mathcal{K}$ is not so far from the usual ERM estimator in
the direct case. Indeed, if $\lambda$ is small enough, we have coarsely:
\[
\tilde h_{G,\lambda}(X_i)=\int_G
\frac{1}{\lambda}\mathcal{K} \biggl(\frac{X_i-x}{\lambda} \biggr)\,\mathrm{d}x
\approx\ind_G(X_i).
\]
As a result, with a small enough bandwidth, the procedure \eqref
{sieve2} reaches the same asymptotic performances as standard ERM.

\section{Conclusion}
\label{section:conclusion}

We have provided in this paper minimax rates of convergence in the
framework of smooth discriminant analysis with errors in variables. In
the presence of plug-in type assumptions, we replace the unknown
densities $f$ and $g$ by deconvolution kernel estimators. It gives a
new family of ERM estimators called deconvolution ERM. It reaches the
minimax rates of convergence. These optimal rates are fast rates
(faster than $n^{-\trup{1}{2}}$) when $\alpha\gamma>d+2\sum_{i=1}^d\beta_i$ and generalize the result of \cite{AT}. As shown in
Table~\ref{tab1}, the influence of the noise $\epsilon$ can be compared with
standard results in nonparametric statistics (see \cite{fan,fan2} for
regression and density estimation with errors in variables or \cite
{Butucea} in goodness-of-fit testing) using kernel deconvolution
estimators. Note that this idea can be adapted to the direct case using
kernel density estimators. It provides a new minimax optimal procedure
in the direct case, under the Plug-in assumption.

\begin{table}[b]
\tablewidth=\textwidth
\tabcolsep=0pt
\caption{Optimal rates of convergence in pointwise $L^2$-risk
in density estimation (see \cite{fan}), optimal separation rates for
goodness-of-fit testing on Sobolev spaces $W(s,L)$ (see, e.g., \cite
{Butucea}) and the result of this work in smooth discriminant analysis
(where $\bar\beta:= \sum_{i=1}^d \beta_i$)}\label{tab1}
\begin{tabular*}{\textwidth}{@{\extracolsep{\fill}}llll@{}}
\hline
&Density estimation& Goodness-of-fit testing &Classification
\\
\hline
Direct case ($\epsilon=0$)& $ n^{-\trup{2\gamma}{(2\gamma+1)}}$ &
$n^{-\trup{2\gamma}{(2\gamma+1/2)}}$ & $n^{-\trup{\gamma(\alpha
+1)}{(\gamma(\alpha+2)+d)}}$
\\
Errors-in-variables & $n^{-\trup{2\gamma}{(2\gamma+2\beta+1)}}$ &
$n^{-\trup{2\gamma}{(2\gamma+2\beta+1/2)}}$ & $n^{-\trup{\gamma
(\alpha+1)}{(\gamma(\alpha+2)+2\bar{\beta}+d)}}$
\\[5pt]
Regularity&$f\in\Sigma(\gamma,L)$& $f\in W(s,L)$ & $f-g\in\Sigma
(\gamma,L)$
\\
\quad assumptions&$|\mathcal{F}[\eta](t)|\sim|t|^{-\beta}$&$|\mathcal
{F}[\eta](t)|\sim|t|^{-\beta}$&$|\mathcal{F}[\eta_i](t)|\sim
|t|^{-\beta_i}$ $\forall i$
\\
\hline
\end{tabular*}
\end{table}

It is important to note that considering the estimation procedure of
this paper, we are facing two different problems of model selection or
adaptation. First of all, the choice of the bandwidths clearly depends
on parameters which may be unknown a priori (e.g., the margin $\alpha$
and the regularity $\gamma$ of the densities). In this sense,
adaptation algorithms should be investigated to choose automatically
$\lambda$ to balance the bias term and the variance term. The second
step of adaptation would be to consider a family of nested $(\mathcal
{G}_k)\subset\mathcal{G}$ and to choose the model which balances the
approximation term and the estimation term. This could be done using
for instance penalization techniques, such as \cite{vdgtsybakov} or
\cite{kolt} or a comparison method such as \cite{tsybakov2004}.

This work can be considered as a first attempt into the study of risk
bounds in classification with errors in variables. It can be extended
in many directions. Naturally the first extension will be to state the
same kind of result in classification. Another natural direction would
be to consider more general complexity assumptions for the hypothesis
space $\mathcal{G}$. In the free-noise case, \cite{localrademacher}
deal with local Rademacher complexities. It allows to consider many
hypothesis spaces, such as VC classes of sets, kernel classes (see
\cite{mendelson}) or even Besov spaces (see \cite{loustau2}). Another
advantage of considering Rademacher complexities is to develop
data-dependent complexities to deal with the problem of model selection
(see \cite{kolt,modelselection}). It also allows us to deal with the
problem of nonunique solution of the empirical minimization.

Into the direction of statistical inverse problem, there are also many
open problems. A natural direction for applications would be to
consider unknown density $\eta$ for the random noise $\epsilon$. This
is a well known issue in the errors-in-variables setting to deal with
unknown operator of inversion. In this setting we can consider repeated
measurements to estimate the density of the noise $\epsilon$ (see, for
instance, \cite{delaiglehallmeister} for both density estimation and
regression with errors). Another natural extension will be to consider
general linear compact operator $A\dvtx f\mapsto Af$ to generalize the case
of deconvolution. In this case, ERM estimators based on standard
regularization methods from the inverse problem literature (see \cite
{engl}) appear as good candidates. This could be the material of future
works.


Finally, the presence of fast rates in discriminant analysis goes back
to \cite{mammen}. In \cite{mammen}, the regularity assumption is
related to the smoothness of the boundaries of the Bayes classifier. If we
consider a set of H\"older boundary fragments, \cite{mammen} states
minimax fast rates in noise-free discriminant analysis. These rates are
attained by ERM estimators. A~natural extension of the present
contribution is to state minimax rates in the presence of H\"older
boundary fragments, where the control of the bias term seems really
more nasty. This is the purpose of a future work.




\section{Proofs}

In this section, with a slight abuse of notations, $C,c,c'>0$ denote
generic constants that may vary from line to line, and even in the same
line. Given two real sequences $(a_n)_{n\in\mathbb{N}}$ and
$(b_n)_{n\in\mathbb{N}}$, the notation $a\simeq b$ (resp. $a\lesssim
b$) means that there exists generic constants $C,c>0$ such that
$ca_n\leq b_n\leq Ca_n$ (resp. $a_n\leq Cb_n$) for all $n\in\mathbb{N}$.\label{section:proofs}

\subsection{Proof of Theorem \texorpdfstring{\protect\ref{thm:lbplugin}}{1}}
\label{s:proof_lb}

The proof mixes standard lower bounds arguments from classification
(see \cite{audibert2004} and \cite{AT}) but then uses some techniques
which are specific to the inverse problem literature (see, for
instance, \cite{Butucea} or \cite{meister}).

Consider $\mathcal{F}_1=\{f_{\rodsigma },\tiksigma =(\sigma_1,\ldots,\sigma_k)\in\{0,1\}
^k\}$ a finite class of densities with respect to a specific measure
$Q_0$ and $g_0$ a fixed density (with respect to the same $Q_0$) such
that $(f_{\rodsigma },g_0)\in\mathcal{F}_{\mathrm
{plug}}$ for all $\tiksigma \in\{0,1\}^k$. The
construction of $f_{\rodsigma }$ as a function of
$\tiksigma $, the value of $g_0$ and the definition of
$Q_0$ will be precised in Section~\ref{section:def}. Then, for all
estimator $\hat G_{n,m}$ of the set $G_K^{\star}$, we have:
%
%
\begin{equation}
\label{eq:borneinf2} \sup_{(f,g)\in\mathcal{F}_{\mathrm{plug}}} \mathbb{E}_{f,g}
d_{\Delta}\bigl(\hat G_{n,m},G_K^{\star}
\bigr) \geq\sup_{f\in\mathcal{F}_1} \mathbb{E}_{g_0} \bigl[
\E_{f} \bigl\{ d_{\Delta}\bigl(\hat G_{n,m},G_K^{\star}
\bigr)|Z^{(2)}_1,\ldots ,Z^{(2)}_m
\bigr\} \bigr].
\end{equation}
In a first time, we propose a triplet $(\mathcal{F}_1,g_0,Q_0)$. Then,
we prove that each associated element satisfies our hypotheses. We
finish the proof with a convenient lower bound for (\ref{eq:borneinf2}).

\subsubsection{Construction of the triplet $(\mathcal{F}_1,g_0,Q_0)$}
\label{section:def}

We only consider the case $d=2$ for simplicity, whereas straightforward
modifications lead to the general $d$-dimensional case. For $g_0$, we
take the constant $1$ over $\R^2$:
\[
g_0(x) = 1, \qquad\forall x\in\R^2.
\]
For any $z\in\R^2$ and positive $\delta$, we write in the sequel
$B(z,\delta):=\{x=(x_1,x_2)\dvtx |x_i-z_{i}|\leq\delta\}$.

For an integer $q\geq1$, introduce the regular grid on $[0,1]^2$
defined as:
\[
G_q= \biggl\{ \biggl(\frac{2p_1+1}{2q},\frac{2p_2+1}{2q}
\biggr),p_i\in\{0,\ldots, q-1\},i=1,2 \biggr\}.
\]
Let $n_q(x)\in G_q$ the closest point to $x\in\R^2$ among points in
$G_q$ (by convention, we choose the closest point to $0$ when it is nonunique). Consider the partition $(\chi_j')_{j=1,\ldots, q^2}$ of
$[0,1]^2$ defined as follows: $x$ and $y$ belongs to the same subset if
and only if $n_q(x)=n_q(y)$. Fix an integer $k\leq q^2$. For any $i\in
\{1,\ldots, k\}$, we define $\chi_i=\chi_i'$ and $\chi_0=\R
^2\backslash\bigcup_{i=1}^k\chi_i$ to get $(\chi_i)_{i=1,\ldots,k}$ a
partition of~$\R^2$. In the sequel, we note by $(z^j)_{j=1,\ldots, k}$
the centers of the $\chi_j$.

Then, we consider the measure $Q_0$ defined as $dQ_0(x)=\mu(x)\,\mathrm{d}x$
where $\mu(x)=\mu_0(x)+\mu_1(x)$ for all $x\in\mathbb{R}^2$ with
\[
\mu_0(x) = k\omega\rho(x_1-1/2)\rho(x_2-1/2)
\quad\mbox{and}\quad\mu_1(x) = (1-k\omega) \rho(x_1-a)
\rho(x_2-b),
\]
where $k$, $\omega$, $a$, $b$ are constants which will be made precise later on and
where for all $x\in\R$, $\rho\dvtx  \R\to[0,1]$ is the function defined as
\[
\rho(x) = \frac{1-\cos(x)}{\pi x^2},\qquad \forall x\in\mathbb{R}.
\]
Recall that $\rho$ satisfies $\mathcal{F}[\rho](t)=(1-|t|)_+$. It
allows us to take advantage of the Noise assumption.
Moreover, $g$ defines a probability density w.r.t. to the measure $Q_0$
since $\int_{\mathbb{R}^2} \mu(x)\,\mathrm{d}x =1$.

Now, we have to define the class $\mathcal{F}_1=\{f_{\rodsigma },\tiksigma \}$. We first introduce $\varphi$
as a $\mathcal{C}^\infty$ probability density function w.r.t. the
measure $Q_0$ and such that
\[
\varphi(x) = 1- c^\star q^{-\gamma}\qquad \forall x
\in[0,1]^2.
\]
Now introduce a class of functions $\psi_j\dvtx \R^2\to\R$, for
$j=1,\ldots, k$ defined for any $x\in\R^2$ as follows:
\[
\psi_j(x)=q^{-\gamma}c_\psi\rho\bigl(2\pi q
\bigl(x_1-z^j_1\bigr)\bigr) \rho\bigl(2\pi q
\bigl(x_2-z^j_2\bigr)\bigr)\cos\bigl(4\pi q
\bigl(x_1-z^j_1\bigr)\bigr)\cos\bigl(4\pi q
\bigl(x_2-z^j_2\bigr)\bigr),
\]
where $(z^j)_{j=1,\ldots, k}$ are the centers of the $\chi_j$.
The class $(\psi_j)_j$ is specific to the noisy case and the inverse
problem literature (see \cite{Butucea} and \cite{meister}). With such
notations, for any $\tiksigma  \in\{0,1 \}
^k$, we define:
\[
f_{\rodsigma }(x)=\varphi(x)+ \sum_{l=1}^k
\sigma _l\psi_l(x), \qquad\forall x\in
\mathbb{R}^2.
\]

Now we have to check that this choice of $\mathcal{F}_1$, $g_0$ and
$Q_0$ provides the Margin assumption and that the complexity assumption
hold true.

\subsubsection{Main assumptions check}

In a first time, we prove that the $f_{\rodsigma }$
define probability density functions w.r.t. the measure $Q_0$. Let
$\tiksigma  \in\{0,1 \}^k$. Remark that,
considering the case $d=1$ w.l.o.g.:
\begin{eqnarray*}
\int_\mathbb{R} \psi_l(x) \mu_0(x)
\,\mathrm{d}x &=& \mathcal{F}[\psi_l \mu _0](0) =
c_\psi q^{-\gamma} \mathcal{F}\bigl[\rho(2\pi q \cdot)
\mu_0(\cdot )\bigr](\pm4\pi q)
\\
&=&c_\psi q^{-\gamma}k\omega\mathcal{F}[\rho]*\mathcal{F}\bigl[
\rho (2\pi q \cdot )\bigr](\pm4\pi q).
\end{eqnarray*}
Then, since
\[
\mathcal{F}\bigl[\rho(2\pi q \cdot )\bigr](t) = \frac{1}{2\pi q} \mathcal{F}[
\rho ] \biggl( \frac{t}{2\pi q} \biggr) \qquad\forall t\in\mathbb{R},
\]
and
\[
\mathcal{F}\bigl[\rho(2\pi q \cdot )\bigr](t) \neq0 \quad\Leftrightarrow\quad-1
< \frac
{t}{2\pi q} <1 \quad\Leftrightarrow\quad-2\pi q < t < 2\pi q,
\]
we get
%
%
\begin{equation}
\label{eq:support} \operatorname{supp} \mathcal{F}[\rho]*\mathcal{F}\bigl[\rho(2\pi
q \cdot )\bigr] = [-2\pi q -1; 2\pi q +1 ] \quad\mbox{and}\quad \int
_\mathbb{R} \psi_l(x) \mu _0(x) \,
\mathrm{d}x =0.
\end{equation}
The same computations show that $\int_\mathbb{R} \psi_l(x) \mu_1(x)
\,\mathrm{d}x =0$ and prove the desired result since $\varphi$ is a probability
density with respect to $Q_0$.

Concerning the regularity, $f_{\rodsigma } \in\Sigma
(\gamma,L)$ for $q$ large enough since $f_{\rodsigma }$
can be written as $q^{-\gamma} F_0(x)$ where $F_0$ is infinitely
differentiable.

In order to conclude this part, we only have to prove that the margin
hypothesis is satisfied for all the couples $(f_{\rodsigma },g)$, namely for some constant $c_2,t_0>0$, we have for $0<t<t_0$:
\[
Q_0 \bigl( \bigl\{ x\in[0,1]^d\dvtx \bigl |
f_{\rodsigma }(x) - g(x) \bigr | \leq t \bigr\} \bigr)\leq c_2
t^{\alpha}.\vadjust{\goodbreak}
\]
First, note that by construction of $Q_0$, we have $\mathrm{d}Q_0(x)=(\mu
_0(x)+\mu_1(x))\,\mathrm{d}x$ and by choosing constant $a,b>0$ large enough in
$\mu_1$, we can restrict ourselves to the study of the Margin
assumption with respect to $Q_0'(\mathrm{d}x)=\mu_0(x)\,\mathrm{d}x$.

Concerning the triplet $(k,\omega,q)$, we set\vspace*{-1pt}
\[
\cases{ %
k = q^2,
\cr
\omega= q^{-\alpha\gamma-2}. }
\]
In particular, we will have $k\omega= q^{-\alpha\gamma}$. Then, we
will distinguish two different cases concerning the possible value of
$t$. The first case concerns the situation where $ C_1 q^{-\gamma} < t
< t_0$ for some constant $C_1$. Then, we have for $Q_0'(\mathrm{d}x)=\mu_0(x)\,\mathrm{d}x$:\vspace*{-1pt}
\[
Q_0' \bigl( \bigl\{ x\in[0,1]^2\dvtx\bigl  |
f_{\rodsigma }(x) - g(x) \bigr | \leq t \bigr\} \bigr)\leq\int_{[0,1]^2}
\mu _0(x) \,\mathrm{d}x \leq k\omega\leq C q^{-\alpha\gamma} \leq C
t^{\alpha}.
\]
Now, we consider the case where $t < C_1 q^{-\gamma}$. For all $\sigma
\in\{0,1\}^{k}$:\vspace*{-1pt}
%
\begin{eqnarray}
\label{marg} Q_0' \bigl( \bigl\{ x
\in[0,1]^2\dvtx \bigl |(f_{\sigma}-g) (x)\bigr |\leq t \bigr\} \bigr)&=&
\int_{[0,1]^2}k\omega\mathbf{1}_{|(f_{\sigma
}-g)(x)|\leq t}\,\mathrm{d}x
\nonumber
\\
&\leq& k\omega\sum_{j=1}^k\int
_{\chi_j}\mathbf{1}_{|(f_{\sigma
}-g)(x)|\leq t}\,\mathrm{d}x
\\[-1pt]
&\leq&k^2\omega\operatorname{Leb} \bigl\{ x\in\chi_{1}
\dvtx \bigl |(f_\sigma -g) (x)\bigr |\leq t \bigr\},
\nonumber
\end{eqnarray}
where without loss of generality, we suppose that $\sigma_1=1$ and we
denote by $\operatorname{Leb}(A)$ the Lebesgue measure of $A$.

Last step is to control the Lebesgue measure of the set $W_1=\{x\in
\chi_1\dvtx |(f_\sigma-g)(x)|\leq t\}$. Since $f_{\sigma}-g=\sum_{j=1}^k\sigma_j\psi_j-c^\star q^{-\gamma}$,
we have\vspace*{-1pt}
\begin{eqnarray*}
W_1 &=& \Biggl\{ x \in\chi_1\dvtx \Biggl\llvert \sum
_{j=1}^k \sigma_j \psi
_j(x) - c^\star q^{-\gamma} \Biggr\rrvert \leq t
\Biggr\}
\\[-1pt]
&=& \Biggl\{ x \in\chi_1\dvtx \Biggl\llvert
\psi_1(x) - \Biggl( c^\star q^{-\gamma} - \sum
_{j=2}^k \sigma_j
\psi_j(x) \Biggr) \Biggr\rrvert \leq t \Biggr\}.
\end{eqnarray*}
Moreover, note that on the square $\chi_j$:\vspace*{-1pt}
%
\begin{eqnarray}
\label{ouff} \sum_{l\neq j}\sigma_l
\psi_l(x)&\leq&q^{-\gamma}c_\psi\sum
_{l\neq j}\frac{1}{2^4\pi^6q^4}\prod_{i=1}^2
\frac
{1}{|x_i-z_{l,i}|^2}
\nonumber
\\[-1pt]
&\leq&\frac{q^{-\gamma}c_\psi}{2^4\pi^6}\sum_{l\neq j}
\frac
{1}{|l-j|^4}
\\[-1pt]
&\leq&\frac{q^{-\gamma}c_\psi}{2^4\pi^6}\zeta(4)=\frac
{q^{-\gamma}c_\psi\pi^4}{90\times2^4\pi^6}:=c'q^{-\gamma},
\nonumber
\end{eqnarray}
where $c'=\frac{c_\psi}{90\times2^4\pi^2}$. Then, if we note by:
\[
c_\infty=\sup_{x\in\chi_1}\rho\bigl(2\pi q
\bigl(x_1-z^1_1\bigr)\bigr) \rho\bigl(2\pi
q\bigl(x_2-z^1_2\bigr)\bigr)\cos\bigl(4\pi q
\bigl(x_1-z^1_1\bigr)\bigr)\cos\bigl(4\pi q
\bigl(x_2-z^1_2\bigr)\bigr),
\]
we have, for any $x\in\chi_1$:
%
\begin{equation}
\label{existenceracine} \sum_{j=1}^k
\sigma_j \psi_j(x)=\psi_1(x)+ \sum
_{j=2}^k \sigma_j \psi_j(x)
\leq\bigl(c_\psi c_\infty+c'\bigr)
q^{-\gamma}.
\end{equation}
Then, for all $x\in\chi_1$, we can define $z^x$ as
\[
z^x = \arg\min_{z: \psi_1(z)= c^\star q^{-\gamma} - \sum_{j=2}^k
\sigma_j \psi_j(z)} \|x-z\|_2.
\]
Indeed, inequality \eqref{existenceracine} ensures the existence of
$z^x$ provided that $c^\star< c_\psi c_\infty+c'$.

In order to evaluate the Lebesgue measure of $W_1$, the main idea is to
approximate $\psi_1$ at each $x\in W_1$ by a Taylor polynomial of
order 1 at $z^x$. We obtain
\begin{eqnarray*}
W_1 & = & \bigl\{ x\in\chi_1\dvtx \bigl |
\psi_1(x) - \psi_1\bigl(z^x\bigr) \bigr |\leq t
\bigr\}
\\
& = & \bigl\{ x\in\chi_1\dvtx \bigl | \bigl\langle D\psi_1
\bigl(z^x\bigr),x-z^x\bigr\rangle +
\psi_1(x) - \psi_1\bigl(z^x\bigr) - \bigl
\langle D\psi_1 \bigl(z^x\bigr),x-z^x\bigr
\rangle\bigr | \leq t \bigr\}
\\
& \subset& \bigl\{ x\in\chi_1\dvtx \bigl\llvert \bigl |\bigl\langle D
\psi_1 \bigl(z^x\bigr),x-z^x\bigr\rangle\bigr | -
\bigl | \psi_1(x) - \psi_1\bigl(z^x\bigr) -\bigl
\langle D\psi_1 \bigl(z^x\bigr),x-z^x\bigr
\rangle\bigr | \bigr\rrvert \leq t \bigr\}.
\end{eqnarray*}
Now, it is possible to see that there exists $c_0>0$ such that
%
\begin{equation}
\label{lowb} \bigl |\bigl\langle D\psi_1 \bigl(z^x
\bigr),x-z^x\bigr\rangle\bigr | \geq c_0 q q^{-\gamma}
\bigl \|x-z^x\bigr \| _1, \qquad \forall x\in\chi_1.
\end{equation}
Moreover, using again the inequality $\no x-z^x\no_1\leq C/q$, there
exists a function $h\dvtx \R\to\R_+$ such that $qh(q)\to0$ as $q\to
\infty$ and which satisfies:
%
\begin{equation}
\label{upb} \frac{| \psi_1(x) - \psi_1(z^x) - \langle D\psi_1
(z^x),x-z^x\rangle|}{\|x-z^x\|_1} \leq q^{-\gamma} h(q).
\end{equation}
At this step, it is important to note that provided that $q:=q(n)\to
\infty$ as $n\to\infty$, there exists some $n_0\in\mathbb{N}$ such
that for any $n\geq n_0$, we have:
\[
\bigl | \bigl\langle D\psi_1 \bigl(z^x\bigr),x-z^x
\bigr\rangle\bigr | > \bigl | \psi_1(x) - \psi_1
\bigl(z^x\bigr) - \bigl\langle D\psi_1
\bigl(z^x\bigr),x-z^x\bigr\rangle\bigr |.
\]
%
Hence, we get the following inclusion
\[
W_1 \subset \biggl\{ x\in\chi_1\dvtx c_0
q q^{-\gamma} \bigl \|x-z^x\bigr \|_1 \biggl(1-
\frac{h(q)}{q} \biggr)\leq t \biggr\}, \qquad\mbox{as } q\rightarrow+\infty.
\]
With the property $qh(q)\to0$ as $q\to\infty$ (or equivalently when
$n\to\infty$), we can find $n_0'$ large enough such that for any
$n\geq n_0'$:
\[
\operatorname{Leb}(W_1) \leq\operatorname{Leb} \biggl( \biggl\{ x
\in\chi _1\dvtx \bigl \|x-z^x\bigr \|_1 \leq
\frac{ t}{2c_0} q^{\gamma-1} \biggr\} \biggr) \leq\frac{t}{2c_0q q^{1-\gamma}}.
\]
Gathering with \eqref{marg}, we hence get, for $t<C_1q^{-\gamma}$,
provided that $\alpha\leq1$:
\begin{eqnarray*}
Q_0' \bigl\{ x\in[0,1]^2\dvtx
\bigl |(f_{\sigma}-g) (x)\bigr |\leq t \bigr\} &\leq& Ck^2\omega
\frac{t}{q^2q^{-\gamma}}
\\
&\leq&Ck\omega\frac{t}{q^{-\gamma}}=Cq^{\gamma(1-\alpha)}t^\alpha
t^{1-\alpha}\leq Ct^\alpha,
\end{eqnarray*}
where $C>0$ is a generic constant.

\subsubsection{Final minoration}
Suppose without loss of generality that $n\leq m$. Now we argue as in
\cite{audibert2004} (Assouad Lemma for classification) and introduce
$\nu$, the distribution of a Bernoulli variable ($\nu(\sigma=1)=\nu
(\sigma=0)=1/2$). Then, denoting by $\P_{\rodsigma }^{\otimes n}$ the law of $(Z^{(1)}_1,\ldots,Z^{(1)}_n)$ when
$f=f_{\rodsigma }$, we get
%
%
\begin{eqnarray}
\label{assouad} && \sup_{\rodsigma \in\{0,1\}}\E_{f} \bigl\{
d_{\Delta
}\bigl(\hat G_{n,m},G_K^{*}
\bigr)|Z^{(2)}_1,\ldots,Z^{(2)}_m \bigr
\}
\nonumber
\\
&&\quad \geq \E_{\nu^{\otimes k}}\mathbb{E}_{f_{\rodsigmas }} d_{\Delta}
\bigl(\hat G_{n,m},G_K^{*}\bigr)
\nonumber
\\
&&\quad \geq \E_{\nu^{\otimes k}}\mathbb{E}_{f_{\rodsigmas }}\sum
_{j=1}^k \int_{\chi_j}\mathbf{1}
\bigl(x\in\hat G_{n,m}\Delta G_K^\star
\bigr)Q_0(\mathrm{d}x)
\\
&&\quad = \sum_{j=1}^k
\E_{\nu^{\otimes(k-1)}}\int_\Omega\E_{\nu
(d\sigma_j)}\int
_{\chi_j}\mathbf{1}\bigl(x\in\hat G_{n,m}(\omega )
\Delta G_K^\star\bigr)Q_0(\mathrm{d}x)
\P_{\rodsigma }^{\otimes
n}(\mathrm{d}\omega)
\nonumber
\\
&&\quad \geq \sum_{j=1}^k
\E_{\nu^{\otimes(k-1)}} \int_\Omega \E_{\nu(d\sigma_j)}\int
_{\chi_j}\mathbf{1}\bigl(x\in\hat G_{n,m}(\omega)\Delta
G_K^\star\bigr)Q_0(\mathrm{d}x) \biggl[
\frac{\P
_{\rodsigma _{j,1}}^{\otimes n}}{\P_{\rodsigma _j}^{\otimes n}} \wedge\frac{\P_{\rodsigma _{j,0}}^{\otimes n}}{\P_{\rodsigma _j}^{\otimes
n}} \biggr]\P_{\rodsigma }^{\otimes n}(
\mathrm{d}\omega),
\nonumber
\end{eqnarray}
where $\tiksigma _{j,r}=(\sigma_1,\ldots,\sigma
_{j-1},r,\sigma_{j+1},\ldots,\sigma_k)$ for $r\in\{0,1\}$.

Moreover, note that from \eqref{ouff}, we have on the square $\chi_j$:
\[
\sum_{l\neq j}\sigma_l
\psi_l(x)\leq c'q^{-\gamma},
\]
where $c'=\frac{c_\psi}{90\times2^4\pi^2}$. Now it is easy to see
that from the definition of the test functions $\psi_j$, for any
integer $k_0,k_1\dvtx k_1>2k_0$, on the square ring $B_j(k_0,k_1)=\{x\in
\chi_j\dvtx \forall i\ |x_i-z_{j,i}|\leq\frac{1}{2k_0q}\mbox{ and }|x_i-z_{j,i}|\geq\frac{1}{k_1q}\}$:
\[
\psi_j(x)-c^\star q^{-\gamma}\geq q^{-\gamma}
\biggl[c_\psi k_0^4\frac{ (1-\cos\trup{2\pi}{k_1} )^2}{\pi^6}(\cos 4
\pi/k_0)^2-c^\star \biggr]= q^{-\gamma},
\]
provided that $c_\psi= \frac{\pi^6(1+c^\star)}{k_0^2(1-\cos2\pi
/k_1)^2(\cos4\pi/k_0)^2}$. Hence, since $c'=\frac{c_\psi}{90\times
2^4\pi^2}$, we can choose $k_0,k_1\in\mathbb{N}$ such that $c'\leq
1$ to get on $B_j(k_0,k_1)$:
%
\begin{equation}
\label{oufff} \sum_{l\neq j}\sigma_l
\psi_l(x)\leq c'q^{-\gamma}\leq\psi
_j(x)-c^\star q^{-\gamma}.
\end{equation}
Now introduce binary valued functions:
\[
\hat{f}(x)=\mathbf{1} (x\in\hat{G}_{n,m})\quad\mbox{and}\quad
f^\star _{\rodsigma }(x)=\mathbf{1} \bigl(x\in G_{K,\sigma}^\star
\bigr),
\]
where $G^\star_{K,\sigma}=\{f_{\rodsigma }-g\geq0\}$.
From \eqref{oufff}, we claim that for any $\tiksigma $:
%
\begin{equation}
\label{fixedsign} \forall x\in B_j(k_0,k_1),
\qquad f^\star_{\rodsigma }(x)= \sigma_j.
\end{equation}
Indeed, since $f_{\rodsigma }-g=\sum_{l=1}^k\sigma
_l\psi_l-c^\star q^{-\gamma}$, gathering with \eqref{oufff}, we have
the following assertion:
\[
f^\star_{\rodsigma }(x)=1\quad\Rightarrow\quad(1+
\sigma_j)\psi _j(x)\geq2c^\star q^{-\gamma}
\quad\Rightarrow\quad\sigma_j=1,
\]
provided that $c^\star\leq q^\gamma\min_{x\in B_j(k_0,k_1)}\psi
_j(x)/2$. Moreover, this choice of $c^\star$ leads to the following assertion:
\[
f^\star_{\rodsigma }(x)=0\quad\Rightarrow\quad\sum
_{l=1}^k\sigma_l\psi_l(x)
\leq c^\star q^{-\gamma}\leq\min_{x\in
B_j(k_0,k_1)}
\psi_j(x)/2.
\]
In this case, if $\sigma_j=1$, we obtain:
%
\begin{equation}
\label{contradiction} \psi_j(x)+\sum_{l\neq j}
\sigma_l\psi_l(x)\leq\min_{x\in
B_j(k_0,k_1)}
\psi_j(x)/2.
\end{equation}
Last step is to show that \eqref{contradiction} is a contradiction.
For this purpose, note that:
\begin{eqnarray*}
\min_{x\in B_j(k_0,k_1)} \biggl(\psi_j(x)+\sum
_{l\neq j}\sigma _l\psi_l(x) \biggr)
&\geq&\min_{x\in B_j(k_0,k_1)}\psi_j(x)+\min_{x\in B_j(k_0,k_1)}
\sum_{l\neq j}\sigma_l\psi_l(x)
\\
&\geq&\min_{x\in B_j(k_0,k_1)}\psi_j(x)/2,
\end{eqnarray*}
where the last inequality is guaranteed when:
\[
\min_{x\in B_j(k_0,k_1)}\psi_j(x)/2\geq-\min
_{x\in
B_j(k_0,k_1)}\sum_{l\neq j}
\sigma_l\psi_l(x).
\]
Finally, the last inequality holds thanks to the positivity of $\psi
_j(x)$ on the set $B_j(k_0,k_1)$ and the fact that $\forall j'\neq j$,
$ \operatorname{sign}\psi_j=\operatorname{sign}\psi_{j'}$. Indeed, $\forall j$,
$\psi_j(x)=0$ for $x\in\mathcal{Z}_{j,1}\cup\mathcal{Z}_{j,2}$ where:
\[
\mathcal{Z}_{j,1}= \biggl\{x\in\R^2\dvtx
\bigl |x_u-z_u^j\bigr |=\frac{l}{q}, u\in\{
1,2\}, l\in\mathbb{N}^* \biggr\}
\]
and
\[
\mathcal{Z}_{j,2}= \biggl\{x\in\R^2\dvtx
\bigl |x_u-z_u^j\bigr |=\frac{2l+1}{8q}, u\in\{1,2
\}, l\in\mathbb{N} \biggr\}.
\]
Note that by construction, $\forall j\neq j'$, $\mathcal
{Z}_{j,2}=\mathcal{Z}_{j',2}=\mathcal{Z}_2$ does not depend on $j\in
\{1,\ldots, k\}$. Moreover, for any $j\in\{1,\ldots, k\}$, $\psi_j$
is alternatively positive and negative on the checkerboard associated
with $\mathcal{Z}_2$. It leads to $\operatorname{sign}\psi_j=\operatorname{sign}\psi_{j'}$, $\forall j\neq j'$ since two centers $z^j$ and
$z^{j'}$ are separated by an odd number of squares (exactly $5$) on
both directions. We hence have by construction that \eqref
{contradiction} is a contradiction and then, \eqref{fixedsign} is shown.

Now we go back to the lower bound. We can write:
\begin{eqnarray*}
\E_{\nu(d\sigma_j)}\int_{\chi_j}\mathbf{1}\bigl(x\in\hat
G_{n,m}(\omega)\Delta G_K^\star
\bigr)Q_0(\mathrm{d}x)&=&\E_{\nu(d\sigma_j)}\int_{\chi_j}
\mathbf{1}\bigl(\hat{f}\neq f^\star_{\rodsigma }\bigr)Q_0(
\mathrm{d}x)
\\
&\geq& \mathbb{E}_{\nu(d\sigma_j)} \biggl[\int_{B_j}\mathbf{1}(\hat{f}
\neq \sigma_j)Q_0(\mathrm{d}x) \biggr]
\\
&=& \frac{1}{2} \biggl[\int_{B_j}\bigl[\mathbf{1}(
\hat{f}\neq1)+\mathbf {1}(\hat{f}\neq0)\bigr]Q_0(\mathrm{d}x)
\biggr]
\\
&=& \frac{1}{2}\int_{B_j}Q_0(x)\,
\mathrm{d}x,
\end{eqnarray*}
where we use \eqref{fixedsign} at the second line with
$B_j:=B_j(k_0,k_1)$. Then it follows from \eqref{assouad} that:
%
%
\begin{eqnarray}
&&{\sup_{\rodsigma \in\{0,+1\}^k}\E_{f} \bigl\{ d_{\Delta}
\bigl(\hat G_{n,m},G_K^{\star}\bigr)|Z^{(2)}_1,
\ldots ,Z^{(2)}_m \bigr\}}
\nonumber
\\
&&\quad \geq \E_{\nu^{\otimes(k-1)}}\sum_{j=1}^k
\int_\Omega \biggl[\frac{\P_{\rodsigma _{j,0}}^{\otimes n}}{\P
_{\rodsigma _j}^{\otimes n}}\wedge\frac{\P
_{\rodsigma _{j,1}}^{\otimes n}}{\P_{\rodsigma _j}^{\otimes n}}
\biggr](\mathrm{d}\omega) \frac{1}{2}\int_{\chi
_j}Q_0(
\mathrm{d}x)\P_{\rodsigma }^{\otimes n}(\mathrm{d}\omega)
\nonumber
\\
&&\quad = \sum_{j=1}^k
\E_{\nu^{\otimes(k-1)}}\bigl[1-\mathbb{V}\bigl(\P _{\rodsigma ,1}^{\otimes n},
\P_{\rodsigma ,0}^{\otimes n}\bigr)\bigr]\frac{1}{2}\int
_{B_j}Q_0(\mathrm{d}x)
\\
&&\quad \geq \sum_{j=1}^k
\E_{\nu^{\otimes(k-1)}}\Bigl[1-\sqrt{\chi^2\bigl(\P
_{\rodsigma ,1}^{\otimes n},\P_{\rodsigma ,0}^{\otimes n}\bigr)}\,\Bigr]
\frac{1}{2}\int_{B_j}Q_0(\mathrm{d}x)
\nonumber
\\
&& \quad= \sum_{j=1}^k\biggl[\biggl(1-\sqrt{
\bigl(1+\chi^2(\P_{1},\P_{0})
\bigr)^n-1}\,\biggr)\frac
{1}{2}\biggr]\int_{B_j}Q_0(
\mathrm{d}x),
\nonumber
\end{eqnarray}
where $\P_{i}$, $i\in\{0,1\}$ is the law of $Z^{(1)}$ when
$f=f_{\rodsigma }$ with $\tiksigma =(i,1,\ldots, 1)$, $i\in\{0,1\}$, $\mathbb{V}(P,Q)$ is the total
variation distance between distribution $P$ and $Q$ and $\chi^2(P,Q)$
is the $\chi^2$ divergence between $P$ and $Q$.
Then we can write, if $\chi^2(\P_{1},\P_{0})\leq\frac{C}{n}$:
%
%
\begin{equation}
\label{fin3} \sup_{\rodsigma \in\{0,+1\}^k} \mathbb {E}_{f_{\rodsigmas },g_0}
d_{\Delta}\bigl(\hat G_{n,m},G_K^{\star}
\bigr)\geq c'\sum_{j=1}^k
\int_{B_j}Q_0(\mathrm{d}x)=c'k
\omega,
\end{equation}
where we use the definition of $Q_0$.

Next step is to find a satisfying upper bound for $\chi^2(\P_{1},\P
_{0})$. We have, by construction of $f_{\rodsigma }$:
\begin{eqnarray*}
\chi^2(\P_{1},\P_0) & = & \int
\frac{ [(f_{\rodsigma ,1}-f_{\rodsigma ,0})\mu*\eta]^2}{f_{\rodsigma ,0}*\eta}\,\mathrm{d}x
\\
& \leq& \int\frac{ [(f_{\rodsigma ,1}-f_{\rodsigma ,0})\mu_0*\eta
]^2}{f_{\rodsigma ,0}\mu*\eta}\,\mathrm{d}x +\int\frac{
[(f_{\rodsigma ,1}-f_{\rodsigma ,0})\mu
_1*\eta]^2}{f_{\rodsigma ,0}\mu*\eta}\,
\mathrm{d}x.
\end{eqnarray*}
The right-hand side term can be considered as negligible with a good
choice of the parameters $a$ and $b$. Hence, we concentrate on the
first one. First, remark that for all $x \in\mathbb{R}^2$, for some $C>0$:
\begin{eqnarray*}
f_{\rodsigma ,0}\mu*\eta&\geq&\frac
{C}{(1+x_1^2)(1+x_2^2)},\qquad \forall x\in
\mathbb{R}^2, \quad\mbox{and}
\\
 \bigl\{(f_{\rodsigma ,+1}-f_{\rodsigma ,0})
\mu_0 \bigr\}*\eta&=& q^{-\gamma} k\omega\{ \psi_l
\rho\}* \eta(x).
\end{eqnarray*}
Then,
\begin{eqnarray*}
\chi^2(\P_1,\P_0) & = & \int
_{\mathbb{R}}\int_{\mathbb{R}} \frac{\{(f_{\omega
_{11}}-f_{\omega_{10}})*\eta(x) \}^2}{f_{\omega_{11}}*\eta
(x)} \,
\mathrm{d}x
\\[1pt]
& \leq& Cq^{-2\gamma}k\omega\int_{\mathbb{R}}\int
_{\mathbb{R}} \bigl(1+x_1^2\bigr)
\bigl(1+x_2^2\bigr) \bigl\{{\psi_1\rho}*
\eta(x) \bigr\}^2 \,\mathrm{d}x.
\end{eqnarray*}
Hence:
\begin{eqnarray*}
\chi^2(\P_1,\P_0) & \leq&
Cq^{-2\gamma}k\omega\int_{\mathbb{R}}\int_{\mathbb{R}}
\bigl\{{\psi_1\rho}*\eta(x) \bigr\}^2 \,\mathrm{d}x +
Cq^{-2\gamma}k\omega \int_{\mathbb{R}}\int_{\mathbb{R}}
x_2^2 \bigl\{{\psi_1\rho }*\eta(x) \bigr
\}^2 \,\mathrm{d}x
\\[1pt]
&&{} + Cq^{-2\gamma}k\omega\int_{\mathbb{R}}\int
_{\mathbb{R}} x_1^2 \bigl\{{
\psi_1\rho}*\eta(x) \bigr\}^2 \,\mathrm{d}x
\\[1pt]
&&{}+
Cq^{-2\gamma
}k\omega\int_{\mathbb{R}}\int_{\mathbb{R}}
x_1^2 x_2^2 \bigl\{ {
\psi_1\rho}*\eta(x) \bigr\}^2 \,\mathrm{d}x
\\[1pt]
& := & A_1 + A_2 + A_3 + A_4.
\end{eqnarray*}
In the following, we only consider the bound of $A_1= C k\omega
q^{-2\gamma} \no(\psi_1 \rho)*\eta\no^2$, the other terms being
controlled in the same way. From the definition of $\psi_1$ and the
conditions on $\eta$, we get
\begin{eqnarray*}
\bigl \no(\psi_1 \rho)*\eta\bigr \no^2&=&\int(\psi_1
\rho)*\eta (x)^2\,\mathrm{d}x=\prod_{i=1}^2
\int\bigl\llvert \mathcal{F}[\psi_1\rho](t_i)\bigr
\rrvert ^2\bigl\llvert \mathcal{F}[\eta_i](t_i)
\bigr\rrvert ^2\,\mathrm{d}t_i
\\
&=& \prod_{i=1}^2\int\bigl\llvert
\mathcal{F}\bigl[\rho(2\pi q\cdot)\rho \bigr](t_i-4\pi q)\bigr
\rrvert ^2\bigl\llvert \mathcal{F}[\eta_i](t_i)
\bigr\rrvert ^2\,\mathrm{d}t_i.
\end{eqnarray*}
Using (\ref{eq:support}), the Noise assumption, and the fact that
$q\rightarrow+\infty$, we get
\begin{eqnarray*}
\bigl \no(\psi_1 \rho)*\eta\bigr \no^2 &=& C q^{-2(\beta_1+\beta_2)} \prod
_{i=1}^2\int\bigl\llvert \mathcal{F}\bigl[
\rho(2\pi q\cdot)\rho\bigr](t_i-4\pi q)\bigr\rrvert
^2\,\mathrm{d}t_i
\\
& = & C q^{-2(\beta_1+\beta_2)} \bigl \|\rho(2\pi q\cdot) \rho\bigr \|^2
\\
& \leq& C q^{-2(\beta_1+\beta_2)} \bigl \|\rho(2\pi q\cdot) \bigr \|^2 \leq C
q^{-2(\beta_1+\beta_2) -2}.
\end{eqnarray*}

Similar bounds are available for $A_2$, $A_3$ and $A_4$ as follows. First,
note that for all $t\in\R$:
\[
\mathcal{F}[\psi_{1}\rho](t) = c_\psi q^{-\gamma}
\mathcal{F}\bigl[\rho (2\pi q\cdot)\rho(\cdot)\bigr](t\pm4\pi q),
\]
and
\[
\frac{\mathrm{d}}{\mathrm{d}t} \mathcal{F}[\psi_{1}\rho](t) = - \bigl(\mathrm{i}
c_\psi q^{-\gamma
}\bigr)^2 t\cdot\mathcal{F}\bigl[
\rho(2\pi q\cdot)\rho(\cdot)\bigr](t\pm4\pi q),
\]
for all $t$ in a subset of $\R$ having a Lebesgue measure equal to $1$.
Then since $\mathcal{F}[\rho]$ and its weak derivative are bounded by
$1$ and supported on $[-1;1]$, we have for instance for $A_2$:
\begin{eqnarray*}
A_2&=&Cq^{-2\gamma}k\omega\int_{\mathbb{R}}\int
_{\mathbb{R}} x_2^2 \bigl\{{
\psi_1\rho}*\eta(x) \bigr\}^2 \,\mathrm{d}x
\\
&\leq& Cq^{-2\gamma}k\omega\int_{\mathbb{R}}\int
_{\mathbb{R}} \biggl( \frac{\mathrm{d}}{\mathrm{d}x_2}{ \mathcal{F}[
\psi_1\rho]}(x)\mathcal {F}[\eta](x) \biggr)^2 \,
\mathrm{d}x,
\end{eqnarray*}
which leads to the same asymptotics as in $A_1$.
It leads to the following upper bound
in the general $d$-dimensional case:
%
%
\begin{equation}
\label{controlchi2} \chi^2(\P_{1},\P_{0})\leq
Cq^{-2\gamma-\alpha\gamma-d-2(\beta
_1+\beta_2)} \leq\frac{C}{n}, \qquad\mbox{with } q= n^{\trup{1}{(2\gamma
+\alpha\gamma+d+2(\beta_1+\beta_2))}}.
\end{equation}
Now using (\ref{fin3}),
\[
\sup_{\sigma\in\{0,1\}^k} \mathbb{E}_{f_{\rodsigmas }} d_{\Delta}\bigl(
\hat G_{n,m},G_K^{\star}\bigr)\geq
c'k\omega=c'q^{-\alpha
\gamma}=c'n^{\trup{-\alpha\gamma}{(2\gamma+\alpha\gamma+d+2(\beta
_1+\beta_2))}},
\]
which concludes the proof of the lower bound.

\subsection{Proof of Theorem \texorpdfstring{\protect\ref{thm:upplugin}}{2}}

The proof is presented for $d=2$ for simplicity whereas straightforward
modifications lead to the $d$-dimensional case. In the sequel, we
identify each $\nu\in\Sigma(\gamma,L)$ with a set $G_\nu= \{
x\dvtx  \nu(x)\geq0 \}$. By the same way, we identify $G_K^\star$
with $\nu^\star= f-g$. Moreover, we assume for simplicity that $n\leq
m$.

\subsubsection{A first inequality}

For all $G_\nu:=\{\nu\geq0\}$, we have, using the notations of
Section~\ref{section:upperbounds}:
\begin{eqnarray*}
&&{R_{n,m}^\lambda(G_\nu) - R_{n,m}^\lambda
\bigl(G_{K}^{\star}\bigr) - R_K^\lambda(G_\nu)+R_K^\lambda
\bigl(G^{\star}_{K}\bigr)}
\\
&& \quad = \frac{1}{2n} \sum_{i=1}^n
U_i(G_\nu) + \frac{1}{2m} \sum
_{i=1}^m V_i(G_\nu) :=
\frac{1}{2} T_{n,m}(G),
\end{eqnarray*}
where, for all $i\in\{1,\ldots, n \}$ and $j\in\{
1,\ldots, m \}$,
\[
U_i(G_\nu) = \bigl\{ h_{K/G_{K}^{\star},\lambda}
\bigl(Z_i^{(1)}\bigr) - h_{G^C_\nu,\lambda}
\bigl(Z_i^{(1)}\bigr) \bigr\}- \mathbb{E} \bigl[
h_{K/G_{K}^{\star},\lambda}\bigl(Z_i^{(1)}\bigr) - h_{G^C_\nu,\lambda
}
\bigl(Z_i^{(1)}\bigr) \bigr],
\]
and
\[
V_j(G_\nu) = \bigl\{ h_{G_{K}^{\star},\lambda}
\bigl(Z_j^{(2)}\bigr) - h_{G_\nu
,\lambda}
\bigl(Z_j^{(2)}\bigr) \bigr\}- \mathbb{E} \bigl[
h_{G_{K}^{\star},\lambda
}\bigl(Z_j^{(2)}\bigr) - h_{G_\nu,\lambda}
\bigl(Z_j^{(2)}\bigr) \bigr].
\]
Then, for all $i\in\{1,\ldots, n \}$, using successively
Lemma~\ref{lip} in the \hyperref[app]{Appendix} and the Margin assumption (Lemma~2 in
\cite{mammen}) we get:
\[
\mathbb{E}\bigl[U_i(G_\nu)\bigr]^2 \leq c
\lambda_1^{-2\beta_1}\lambda _2^{-2\beta_2}
d_{\Delta}\bigl(G_\nu,G_{K}^{\star}\bigr)
\leq c'\lambda _1^{-2\beta_1}\lambda_2^{-2\beta_2}
d_{f,g}\bigl(G_\nu,G_{K}^{\star
}
\bigr)^{\trup{\alpha}{(\alpha+1)}},
\]
and
\[
\bigl | U_i(G_\nu) \bigr | \leq C \prod
_{i=1}^2 \lambda_i^{-\beta_i-1/2},
\]
for some constant $C>0$. The Bernstein's inequality leads to
\begin{eqnarray*}
&&P \Biggl( \Biggl\llvert \frac{1}{n} \sum_{i=1}^n
U_i(G_\nu) \Biggr\rrvert >a \Biggr)
\\
&&\quad\leq2 \exp \biggl[
- \frac{Cna^2}{ a \times\lambda
_1^{-\beta_1-1/2}\lambda_2^{-\beta_2-1/2} + \lambda_1^{-2\beta
_1}\lambda_2^{-2\beta_2} d_{f,g}(G_\nu,G_{K}^{\star})^{\trup
{\alpha}{(\alpha+1)}}} \biggr],
\end{eqnarray*}
for all $a>0$. Since $\beta_i>1/2$ for all $i\in\{1,\ldots,d
\}$, the particular choice $a=d_{f,g}(G_\nu,G_{K}^{\star})$ yields
\begin{eqnarray*}
P \Biggl( \Biggl\llvert \frac{1}{n} \sum_{i=1}^n
U_i(G_\nu) \Biggr\rrvert >d_{f,g}
\bigl(G_\nu,G_K^\star\bigr) \Biggr) & \leq& 2
\exp \bigl[ - Cn \lambda_1^{2\beta_1}\lambda_2^{2\beta
_2}
d_{f,g}\bigl(G_\nu,G_K^\star
\bigr)^{2- \trup{\alpha}{(\alpha+1)}} \bigr]
\\
& = & 2 \exp \bigl[ - Cn \lambda_1^{2\beta_1}
\lambda_2^{2\beta_2} d_{f,g}\bigl(G_\nu,G_K^\star
\bigr)^{\trup{(2+\alpha)}{(\alpha+1)}} \bigr].
\end{eqnarray*}
In the upper bound above, we have implicitly use the fact that
$d_{f,g}(G_\nu,G_K^\star)/2 \leq(d_{f,g}(G_\nu,G_K^\star
)/\allowbreak  2)^{\alpha/(\alpha+1)}$ since $d_{f,g}(G_1,G_2)\leq2$ for all
$G_1,G_2 \subset K$. Using the same algebra on the $V_j(G_\nu)$, we get
\[
P \bigl( \bigl\llvert T_{n,m}(G_\nu) \bigr\rrvert
>d_{f,g}\bigl(G_\nu,G_K^\star\bigr)
\bigr) \leq2 \exp \bigl[ - Cn \lambda_1^{2\beta_1}\lambda
_2^{2\beta_2} d_{f,g}\bigl(G_\nu,G_K^\star
\bigr)^{\trup{(2+\alpha)}{(\alpha
+1)}} \bigr].
\]
This concludes the first part of the proof. Let $t$ a positive
parameter which will be chosen further and introduce the set $\mathcal
{G}'$ defined as
\[
\mathcal{G}'= \bigl\{ G\in\mathcal{N}_{\delta_n},
d_{f,g}\bigl(G_K^\star,G\bigr) > t
\delta_n^{1+\alpha} \bigr\},
\]
where $\mathcal{N}_{\delta_n}$ is the $\delta_n$-network introduced
in Section~\ref{section:upper_bound}, with $\delta_n=\delta_{n,n}$.
Using the upper bound above,
\begin{eqnarray*}
&&P \biggl( \exists G\in\mathcal{G}'\dvtx \bigl |T_{n,m}(G)\bigr |
\geq\frac{1}{4} d_{f,g}\bigl(G,G_K^\star
\bigr) \biggr)
\\
&&\quad \leq  \sum_{G\in\mathcal{G}'} P \biggl(
\bigl |T_{n,m}(G)\bigr | \geq\frac
{1}{4} d_{f,g}
\bigl(G,G_K^\star\bigr) \biggr)
\\
&& \quad \leq  \sum_{G\in\mathcal{G}'} 2 \exp \bigl[ - Cn \lambda
_1^{2\beta_1}\lambda_2^{2\beta_2}
d_{f,g}\bigl(G,G_K^\star\bigr)^{\trup
{(2+\alpha)}{(\alpha+1)}}
\bigr]
\\
&&\quad  \leq  \sum_{G\in\mathcal{G}'} 2 \exp \bigl[ - Cn \lambda
_1^{2\beta_1}\lambda_2^{2\beta_2} \bigl(t
\delta_n^{1+\alpha} \bigr)^{\trup{(2+\alpha)}{(\alpha+1)}} \bigr]
\\
&& \quad \leq  \sum_{G\in\mathcal{G}'} 2 \exp \bigl[ - Cn \lambda
_1^{2\beta_1}\lambda_2^{2\beta_2}
t^{\trup{(2+\alpha)}{(\alpha+1)}} \delta_n^{2+\alpha} \bigr].
\end{eqnarray*}
Since $\log\operatorname{card}(\mathcal{N}_{\delta_n}) \leq A\delta
_n^{-2/\gamma}$, we get
\[
P \bigl( \exists G\in\mathcal{G}'\dvtx \bigl |T_{n,m}(G)\bigr |
\geq\tfrac{1}{4} d_{f,g}\bigl(G,G_K^\star
\bigr) \bigr) \leq\exp \bigl[ A\delta_n^{-2/\gamma
} - Cn
\lambda_1^{2\beta_1}\lambda_2^{2\beta_2}
t^{\trup{(2+\alpha
)}{(\alpha+1)}} \delta_n^{2+\alpha} \bigr].
\]
Thanks to the value of $\delta_n$, we get $ \delta_n^{-2/\gamma}
\simeq n\lambda_1^{2\beta_1}\lambda_2^{2\beta_2} \delta
_n^{2+\alpha}$.
Hence, for $t$ large enough,
\begin{eqnarray}
\label{eq:dernier}
&&P \biggl( \exists G\in\mathcal{G}'\dvtx
\bigl |T_{n,m}(G)\bigr | \geq\frac{1}{4} d_{f,g}
\bigl(G,G_K^\star\bigr) \biggr)
\nonumber
\\
&&\quad\leq  \exp \bigl[ - C t
\delta_n^{-2/\gamma} \bigr]
\\
&&\quad=\exp \biggl[ - Ct \biggl( \frac{\lambda_1^{-\beta_1}\lambda
_2^{-\beta_2}}{ \sqrt{n}} \biggr)^{-(\trup{2}{\gamma})(\trup
{2}{(2/\gamma+2+\alpha)})} \biggr].
\nonumber
\end{eqnarray}
Now, using Lemma~\ref{link} in the \hyperref[app]{Appendix}, we can find a set $G_n \in
\mathcal{N}_{\delta_n}$ such that:
\[
d_{f,g}\bigl(G_n,G_K^\star\bigr)
\leq c_2\bigl \no\nu^*-\nu_n\bigr \no_\infty^{\alpha
+1}
\leq c_2\delta_n^{1+\alpha}.
\]
Then, for all $G \in\mathcal{G}'$, we get
\[
\frac{1}{8} d_{f,g}\bigl(G,G_K^\star
\bigr) - \frac{3}{4} d_{f,g}\bigl(G_n,G_K^\star
\bigr) \geq\frac{t}{8} \delta_n^{1+\alpha} -
\frac{3c_2}{4} \delta _n^{1+\alpha} \geq\frac{c_2}{4}
\delta_n^{1+\alpha},
\]
provided that $t>8c_2$. We eventually obtain:
%
%
\begin{eqnarray}
\label{eq:mino2} && P \bigl( d_{f,g}\bigl(G_{K}^{\star},
\hat G_{n,m}\bigr) > t \delta _n^{1+\alpha} \bigr)
\nonumber
\\
&&\quad\leq P \bigl( \exists G\in\mathcal{G}' \dvtx
R^\lambda_{n,m}(G) \leq R^\lambda_{n,m}(G_n)
\bigr)
\\
&&\quad = P \bigl( \exists G\in\mathcal{G}' \dvtx
\tfrac{1}{2}d_{f,g}^\lambda \bigl(G,G_K^\star
\bigr)+ T_{n,m}(G) - \tfrac{1}{2}d_{f,g}^\lambda
\bigl(G_n,G_K^\star \bigr) -
T_{n,m}(G_n) \leq0 \bigr),
\nonumber
\end{eqnarray}
where for all $G_1, G_2 \subset K$,
\[
\tfrac{1}{2}d_{f,g}^\lambda(G_1,G_2):
= R^\lambda_K(G_1)-R^\lambda_K(G_2).
\]

\subsubsection{Control of the bias}

Last step is to control the bias term. In particular, given $G_1,G_2
\subset K$, we want to measure the difference between
$R_K(G_1)-R_K(G_2)$ and $R^\lambda_K(G_1)-R^\lambda_K(G_2)$. First of
all, we have to explicit the term $R_K^\lambda$. Recall that for all
$G_1 \subset K$,
\begin{eqnarray*}
2R_K^\lambda(G_1) & := & 2 \mathbb{E}
R_{n,m}^\lambda(G_1)
\\
& = & \mathbb{E}\bigl[h_{K/G_1,\lambda} \bigl(Z_1^{(1)}
\bigr)\bigr] + \mathbb {E}\bigl[h_{G_1,\lambda} \bigl(Z_1^{(2)}
\bigr)\bigr]
\\
& = & \mathbb{E} \biggl[\int_{ K/G_1} \frac{1}{\lambda}
\mathcal {K}_\eta \biggl( \frac{Z_1^{(1)}-x}{\lambda} \biggr)\,\mathrm{d}x
\biggr] +\mathbb{E} \biggl[\int_{ G_1} \frac{1}{\lambda}
\mathcal{K}_\eta \biggl( \frac{Z_1^{(2)}-x}{\lambda} \biggr)\,\mathrm{d}x
\biggr]
\\
& = & \int_{ K/G_1} \mathbb{E} \biggl[ \frac{1}{\lambda}
\mathcal {K}_\eta \biggl( \frac{X_1^{(1)}+\epsilon_1^{(1)}-x}{\lambda} \biggr) \biggr]\,
\mathrm{d}x +\int_{ G_1} \mathbb{E} \biggl[ \frac{1}{\lambda}
\mathcal{K}_\eta \biggl( \frac{X_1^{(2)}+\epsilon_1^{(2)}-x}{\lambda
} \biggr) \biggr]\,
\mathrm{d}x.
\end{eqnarray*}
Using the properties of the deconvolution kernel, we can see that for
all $x \in K$,
\[
\mathbb{E} \biggl[ \frac{1}{\lambda} \mathcal{K}_\eta \biggl(
\frac
{X_1^{(1)}+\epsilon_1^{(1)}-x}{\lambda} \biggr) \biggr] = \mathbb {E} \biggl[ \frac{1}{\lambda}
\mathcal{K} \biggl( \frac
{X_1^{(1)}-x}{\lambda} \biggr) \biggr] = \int_{\mathbb{R}^d}
\frac
{1}{\lambda} \mathcal{K} \biggl( \frac{y-x}{\lambda} \biggr) f(y) \,
\mathrm{d}Q(y).
\]
The same result holds true when replacing $X_1^{(1)}$ by $X_1^{(2)}$
and $f$ by $g$. Hence, we obtain that
\[
2R_K^\lambda(G_1) = \int_{ K/G_1}
\int_{\mathbb{R}^d} \frac
{1}{\lambda} \mathcal{K} \biggl(
\frac{y-x}{\lambda} \biggr) f(y) \,\mathrm{d}Q(y) \,\mathrm{d}x +\int
_{ G_1} \int_{\mathbb{R}^d} \frac{1}{\lambda}
\mathcal{K} \biggl( \frac{y-x}{\lambda} \biggr) g(y) \,\mathrm{d}Q(y)\,\mathrm{d}x.
\]
%
Moreover, if $Q$ is not the Lebesgue measure, note that by assumption,
there exists a constant $c_0>0$ such that:
%
\begin{equation}
\label{c0togetd} \int_{G_1\Delta G_2} \mathrm{d}x\leq
c_0^{-1}d_\Delta(G_1,G_2).
\end{equation}
We then have:
\begin{eqnarray*}
&&{\bigl\llvert \bigl(R^\lambda_K-R_K\bigr)
(G_1-G_2) \bigr\rrvert }
\\
&&\quad \leq \frac{1}{2} \biggl\llvert \int \biggl[ \int\frac{1}{\lambda
}
\mathcal{K} \biggl(\frac{y-x}{\lambda} \biggr)f(y)\mu(y)\,\mathrm{d}y-f(x)\mu (x)
\biggr] \bigl[\mathbf{1} (x\in K/G_1)-\mathbf{1} (x\in
K/G_2) \bigr]\,\mathrm{d}x
\\
&&\qquad{} +   \int \biggl[\int\frac{1}{\lambda}\mathcal{K} \biggl(
\frac{y-x}{\lambda} \biggr)g(y)\mu(y)\,\mathrm{d}y-g(x)\mu(x) \biggr] \bigl[
\mathbf{1} (x\in G_1)-\mathbf{1} (x\in G_2) \bigr]\,
\mathrm{d}x \biggr\rrvert
\\
&&\quad \leq \frac{1}{2} \int_{G_1\Delta G_2} \bigl\llvert
\mathcal{K}_\lambda *\bigl(\nu^\star\cdot\mu\bigr) (x)-
\nu^\star\cdot\mu(x)\bigr\rrvert \,\mathrm{d}x
\\
&&\quad \leq \frac{1}{2c_0}\bigl \| \mathcal{K}_\lambda*\bigl(
\nu^\star\cdot\mu\bigr) - \nu^\star\cdot\mu\bigr \|_\infty
\int_{G_1\Delta G_2} \,\mathrm{d}x
\\
&&\quad \leq C d_\Delta(G_1,G_2) \bigl[
\lambda_1^\gamma+\lambda _2^\gamma
\bigr]
\\
&&\quad \leq C \bigl[\lambda_1^\gamma+
\lambda_2^\gamma \bigr]d_{f,g}(G_1,G_2)^{\trup{\alpha}{(\alpha+1)}},
\end{eqnarray*}
for some $C>0$, where $\mathcal{K}_\lambda(\cdot)=\frac{1}{\lambda
} \mathcal{K}(\cdot/\lambda)$. Indeed, provided that $\nu\mu\in
\Sigma(\gamma,L)$ and $\mathcal{K}$ is a kernel of order $l=\lfloor
\gamma\rfloor$, it is well known that:
%
\begin{equation}
\label{holdertrick} \bigl \| \mathcal{K}_\lambda*(\nu\mu) - \nu\mu
\bigr \|_\infty\leq C \bigl[\lambda_1^\gamma+
\lambda_2^\gamma \bigr].
\end{equation}
This bound is sufficient for the case $\alpha=0$. If $\alpha>0$,
using the Young inequality:
\[
xy^r \leq ry + (1-r)x^{1/(1-r)},\qquad \forall x,y \in
\mathbb{R}^+,
\]
with $r=\alpha/(\alpha+1)$, $x=C\kappa^{-\alpha/\alpha+1}
[\lambda_1^\gamma+\lambda_2^\gamma ]$ and $y=\kappa d_{f,g}
(G_1,G_2)$, where $\kappa>0$ is chosen later on, we get for all
$G_1,G_2 \subset K$:
%
%
\begin{equation}
\label{eq:mino3} \bigl\llvert \bigl(R^\lambda_K-R_K
\bigr) (G_1-G_2) \bigr\rrvert \leq \biggl(1-
\frac{\alpha}{\alpha+1} \biggr) \biggl(\frac
{C}{\kappa} \biggr)^{\alpha} \bigl[
\lambda_1^\gamma+\lambda _2^\gamma
\bigr]^{\alpha+1} + \frac{\alpha}{\alpha+1}\kappa d_{f,g}
(G_1,G_2).\quad
\end{equation}

\subsubsection{Conclusion of the proof}

Hence, it follows from (\ref{eq:mino2}) and (\ref{eq:mino3}) that if
$\alpha>0$, by choosing $\kappa=(\alpha+1)/(4\alpha)$:
\begin{eqnarray*}
&& {P \bigl( d_{f,g}\bigl(G_K^\star,\hat
G_{n,m}\bigr) > t \delta _n^{1+\alpha} \bigr)}
\\
&&\quad \leq P \Biggl( \exists G\in\mathcal{G}' \dvtx \biggl(
\frac{1}{2} - \frac{\alpha}{\alpha+1}\kappa \biggr) d_{f,g}
\bigl(G,G_K^\star\bigr) + T_{n,m}(G)
\\
&&\phantom{\quad \leq P \Biggl(} {}   - \biggl( \frac{1}{2} +
\frac{\alpha}{\alpha+1}\kappa \biggr) d_{f,g}\bigl(G_n,G_K^\star
\bigr) - T_{n,m}(G_n) + C\sum_{i=1
}^2
\lambda_i^{\gamma(1+\alpha) } \leq0 \Biggr)
\\
&&\quad \leq P \biggl( \exists G\in\mathcal{G}'\dvtx
T_{n,m}(G) \leq-\frac
{1}{8} d_{f,g}
\bigl(G,G_K^\star\bigr) \biggr) + P \Biggl(
T_{n,m}(G_n) \geq C \Biggl(\delta_n^{1+\alpha}
+ \sum_{i=1 }^2\lambda_i^{\gamma
(1+\alpha) }
\Biggr) \Biggr).
\end{eqnarray*}
Note that the same inequalities hold for $\alpha=0$ using the crude bound:
\[
\bigl\llvert \bigl(R^\lambda_K-R_K\bigr)
(G_1-G_2) \bigr\rrvert \leq C\bigl[\lambda_1^\gamma
+\lambda_2^\gamma\bigr].
\]
In order to conclude, remark that the proposed choice of $(\lambda
_j)_{j=1,2}$ provides:
\[
\delta_n^{1+\alpha} \simeq\sum_{i=1 }^2
\lambda_i^{\gamma(1+\alpha
) } \quad\Leftrightarrow\quad\forall i\in\{1,2
\},\qquad \biggl( \frac{\lambda_1^{-\beta
_1}\lambda_2^{-\beta_2}}{ \sqrt{n}} \biggr)^{\trup{(2\gamma(\alpha
+1))}{(\gamma(2 + \alpha)+2)}} \simeq\lambda_i^{\gamma(\alpha+1)}.
\]
Using (\ref{eq:dernier}), we eventually get
\begin{eqnarray*}
&&{P \bigl( d_{f,g}\bigl(G_K^\star,\hat
G_{n,m}\bigr) > t n^{-\trup{(
\gamma(\alpha+1))}{( \gamma(2+\alpha)+2+ 2\sum_{i=1}^{2}\beta_i )}} \bigr)}
\\
&&\quad \leq \exp \bigl[ - C_1t n^{\trup{1 }{( \gamma(2+\alpha)+2+ 2\sum
_{i=1}^{2}\beta_i )}} \bigr]+ \exp
\bigl[ - C_2 n^{\trup{1 }{(
\gamma(2+\alpha)+2+ 2\sum_{i=1}^{2}\beta_i )}} \bigr],
\end{eqnarray*}
where $C_1$, $C_2$ denote positive constants. In order to conclude, we can
remark that
\begin{eqnarray*}
&& {n^{\tau_d(\alpha,\beta,\gamma)} \mathbb{E}_{f,g} d_{f,g}
\bigl(G_K^\star,\hat G_{n,m}\bigr)}
\\
&& \quad\leq t + \mathbb{E}_{f,g} d_{f,g}
\bigl(G_K^\star,\hat G_{n,m}\bigr) \mathbf
{1}_{ \{ d_{f,g}(G_K^\star,\hat G_{n,m}) > t n^{-\trup{
\gamma(\alpha+1)}{( \gamma(2+\alpha)+d+ 2\sum_{i=1}^{d}\beta_i )}}
\}}
\\
&& \quad\leq t +2\exp \bigl[ - C_1t n^{\trup{1 }{( \gamma(2+\alpha)+2+
2\sum_{i=1}^{2}\beta_i )}} \bigr]+ 2\exp
\bigl[ - C_2 n^{\trup{1
}{( \gamma(2+\alpha)+2+ 2\sum_{i=1}^{2}\beta_i) }} \bigr] \leq C
\end{eqnarray*}
for some positive constant $C$, where we have used the bound
$d_{f,g}(G_1,G_2) \leq2$ for all \mbox{$G_1,G_2 \subset K$}.

\subsection{Proof of Corollary \texorpdfstring{\protect\ref{cor}}{1}}

The proof follows the same steps as the proof of Theorem~\ref
{thm:upplugin}. Note that in the direct case, using a kernel $\mathcal
{K}$ with bounded Fourier transform, we have under the Margin assumption:
\[
\E\bigl[U_i(G)^2\bigr]\leq Cd_\Delta
\bigl(G,G_K^\star\bigr)\leq C'd_{f,g}
\bigl(G,G_K^\star \bigr)^{\trup{\alpha}{(\alpha+1)}}  \quad\mbox{and}
\quad \bigl |U_i(G)\bigr |\leq C,
\]
for some constant $C>0$. Remark that the last inequality is more
precise than in the error-in-variable case. Then using Bernstein's
inequality, we have exactly as in the proof of Theorem~\ref{thm:upplugin}:
\[
P \Biggl( \Biggl\llvert \frac{1}{n} \sum_{i=1}^n
U_i(G_\nu) \Biggr\rrvert >a \Biggr) \leq2 \exp \biggl[
- \frac{Cna^2}{ a + d_{\Delta}(G_\nu
,G_{K}^{*})} \biggr],
\]
for all $a>0$. Choosing $a=d_{f,g}(G_\nu,G_K^\star)$ and using the
same algebra, we get a control of the upper bound provided that:
\[
\delta_n^{-\trup{2}{\gamma}} \simeq n\delta_n^{2+\alpha}
\quad\mbox{and}\quad \delta_n^{1+\alpha} \geq\sum
_{i=1 }^2\lambda_i^{\gamma
(1+\alpha) }.
\]
The choice of $\lambda$ and $\delta_n$ in Corollary~\ref{cor}
concludes the proof.

\begin{appendix}\label{app}
%
%
\section*{Appendix}

\begin{lemma}
\label{link}
For any $(f,g)$ satisfying the Margin assumption with parameter $\alpha
>0$, we have:
\[
d_{f,g}\bigl(G_\nu,G_K^\star\bigr)
\leq c_2 \bigl \no\nu-\nu^\star\bigr \no_\infty
^{\alpha+1},
\]
where $G_\nu=\{\nu\geq0\}$ and $\nu^\star=f-g$.
\end{lemma}

\begin{pf}
The proof is a straightforward modification of the
proof of Lemma~5.1 in \cite{AT} which state a similar result in the
binary classification framework. In the following, given $x\in\mathbb
{R}$, we write $\operatorname{sign}(x)=1$ if $x>0$, $\operatorname{sign}(x)=0$ if
$x=0$, and $\operatorname{sign}(x)=-1$ if $x<0$. Then, we get
\begin{eqnarray*}
d_{f,g}\bigl(G_\nu,G_K^\star\bigr)
& = & \int_K \bigl |\nu^\star(x) \bigr |
\mathbf{1}_{ \{ x\in
G_K^\star\Delta G_\nu \}}\,\mathrm{d}Q(x)
\\
& = & \int_K \bigl |\nu^\star(x)\bigr  |
\mathbf{1}_{ \{\operatorname{sign}(\nu^\star(x)) \neq\operatorname{sign}(\nu(x)) \}
}\,\mathrm{d}Q(x)
\\
& \leq& \int_K \bigl |\nu^\star(x) \bigr |
\mathbf{1}_{ \{0<|\nu
^\star(x)| \leq|\nu(x)-\nu^\star(x)|  \}}\,\mathrm{d}Q(x)
\\
& \leq& \bigl \| \nu- \nu^\star\bigr \|_\infty Q \bigl( \bigl\{ x\in K
\dvtx 0<\bigl |\nu ^\star(x)\bigr | \leq\bigl \|\nu-\nu^\star\bigr \|_\infty
\bigr\} \bigr) \leq c_2 \bigl \no\nu-\nu^\star
\bigr \no_\infty^{\alpha+1},
\end{eqnarray*}
where we have used the Margin assumption in order to get the last inequality.
\end{pf}

\begin{lemma}
\label{lip}
Assume that $\eta$ satisfies the {Noise assumption}. Let
$\mathcal{K}_{\eta}$ a deconvolution kernel defined in \eqref{dk}
such that $ \mathcal{F}[\mathcal{K}]$ is bounded and compactly
supported. If $Q(x)=\mu(x)\,\mathrm{d}x$, we assume that $\min_{x\in K} \mu(x)
\geq c_0$ for some $c_0>0$.
Then, we have,
\[
\begin{array} {@{}r@{\quad}l@{}} (\mathrm{i}) & \E\bigl[h_{G,\lambda}(Z)-h_{G',\lambda}(Z)
\bigr]^2\leq Cd_\Delta \bigl(G,G'\bigr)
\displaystyle\prod_{i=1}^d
\lambda_i^{-2\beta_i},
\\
\noalign{\vspace*{2pt}} (\mathrm{ii}) & \displaystyle\sup_{x\in K}
\bigl | h_{G,\lambda}(x) - h_{G',\lambda}(x) \bigr | \leq C\prod
_{i=1}^d \lambda_i^{-\beta_i-1/2},
\end{array}
\]
for some generic constant $C>0$.
\end{lemma}

\begin{pf}
For the sake of convenience, we only consider the case
where $d=1$. We first prove (i). We have, using \eqref{c0togetd}:
\begin{eqnarray*}
&&\E\bigl[h_{G,\lambda}(Z)-h_{G',\lambda}(Z)\bigr]^2
\\
&&\quad  =  \int
_\mathbb{R} \biggl[ \int_{\mathbb{R}}
\frac{1}{\lambda} \mathcal{K}_\eta \biggl(\frac{z-x}{\lambda} \biggr) (
\mathbf {1}_{\{ x\in G\}}-\mathbf{1}_{\{ x\in G'\}
})\mathbf{1}_{\{ x\in K\}}
\,\mathrm{d}x \biggr]^2 (f\mu)*\eta (z)\,\mathrm{d}z
\\
 &&\quad \leq  c\int_{\mathbb{R}} \frac{1}{\lambda^2} \bigl\llvert \mathcal
{F}\bigl[\mathcal{K}_\eta(\cdot/\lambda)\bigr](t) \bigr\rrvert
^2 \bigl\llvert \mathcal {F}\bigl[ (\mathbf{1}_{\{\cdot\in G\}}-
\mathbf{1}_{\{\cdot\in
G'\}})\mathbf{1}_{\{\cdot\in K\}}\bigr](t) \bigr\rrvert
^2 \,\mathrm{d}t
\\
 &&\quad \leq  C \max_{x\in\mathbb{R}^d} \mu(x) \times\lambda^{-2\beta
}
\int_{K} \mathbf{1}_{\{ t \in G\Delta G'\}} \,\mathrm{d}t
\\
 &&\quad \leq  C\lambda^{-2\beta} d_{\Delta}\bigl(G,G'
\bigr).
\end{eqnarray*}
Indeed, for all $s\in\mathbb{R}$, using assumptions on the kernel
$\mathcal{K}_\eta$:
%
\setcounter{equation}{0}
\begin{eqnarray}
\label{eq:inter}
 \frac{1}{\lambda^2} \bigl\llvert \mathcal{F}\bigl[
\mathcal{K}_\eta(\cdot /\lambda)\bigr](s) \bigr\rrvert ^2 &=&
\bigl\llvert \mathcal{F}[\mathcal{K}_\eta ](s\lambda) \bigr\rrvert
^2 =\biggl\llvert \frac{\mathcal{F}[\mathcal
{K}](s\lambda)}{\mathcal{F}[\eta](s)} \biggr\rrvert ^2
\nonumber
\\[-8pt]
\\[-8pt]
&\leq& C\sup_{s\in
[-M/\lambda,M/\lambda]} \biggl\llvert \frac{1}{\mathcal{F}[\mathcal
{K}_\eta](s)} \biggr\rrvert
^2 \leq C \lambda^{-2\beta},\nonumber
\end{eqnarray}
where $\mathcal{F}[\mathcal{K}]=0$ on $\R\setminus[-M,M]$.

In order to prove  (ii), we use the following algebra
\begin{eqnarray*}
\sup_{z\in\mathbb{R}} \bigl | h_{G,\lambda}(z) - h_{G',\lambda}(z)\bigr  | &
\leq& \sup_{z\in\mathbb{R}} \int_{G\Delta G'}
\frac{1}{\lambda
} \biggl\llvert \mathcal{K}_\eta \biggl(
\frac{z-x}{\lambda} \biggr)\biggr\rrvert \,\mathrm{d}x
\\
& \leq& C\sup_{z\in\mathbb{R}} \int_{K}
\frac{1}{\lambda} \biggl\llvert \mathcal{K}_\eta \biggl(
\frac{z-x}{\lambda} \biggr)\biggr\rrvert \,\mathrm{d}x
\\
& \leq& C \sup_{z\in\mathbb{R}} \sqrt{ \int\frac{1}{\lambda^2}
\mathcal{K}^2_\eta \biggl(\frac{z-x}{\lambda} \biggr) \,
\mathrm{d}x}
\\
&\leq& C \lambda^{-1/2}\sqrt{ \int_{[-M,M]}\biggl
\llvert \frac{\mathcal
{F}[\mathcal{K}](t)}{\mathcal{F}[\eta](t/\lambda)}\biggr\rrvert ^2 \,\mathrm{d}t}
\\
&\leq& C \lambda^{-\beta-1/2},
\end{eqnarray*}
where last line uses the Noise assumption and assumptions on the kernel
$\mathcal{K}_\eta$.
\end{pf}
\end{appendix}

\section*{Acknowledgements}

We would like to thank both referees and the Associate Editor whose
remarks and valuable comments help considerably to improve the paper.


%

\printhistory

\end{document}